\documentclass[10pt,english]{article}
\usepackage[T1]{fontenc}
\usepackage[latin1]{inputenc}
\usepackage{a4wide}
\usepackage{babel}
\usepackage{psfrag}

\makeatletter

\newtheorem{prop}{Proposition}[section]
\newtheorem{lemma}[prop]{Lemma}
\newtheorem{theorem}[prop]{Theorem}
\newtheorem{defi}[prop]{Definition}
\newtheorem{rem}[prop]{\em Remark\/}
\newtheorem{cor}[prop]{Corollary}
\newtheorem{cond}[prop]{Conditions}
\newtheorem{ex}[prop]{\em Example\/}

\newcommand{\bp}{\begin{prop}}
\newcommand{\bl}{\begin{lemma}}
\newcommand{\bt}{\begin{theorem}}
\newcommand{\bd}{\begin{defi}\rm}
\newcommand{\br}{\begin{rem}\rm}
\newcommand{\be}{\begin{equation}}
\newcommand{\bea}{\begin{eqnarray}}
\newcommand{\bpr}{\begin{proof}}
\newcommand{\bc}{\begin{cor}}
\newcommand{\bco}{\begin{cond}}
\newcommand{\bex}{\begin{ex}\rm}

\newcommand{\ep}{\end{prop}}
\newcommand{\el}{\end{lemma}}
\newcommand{\et}{\end{theorem}}
\newcommand{\ed}{\end{defi}}
\newcommand{\er}{\end{rem}}
\newcommand{\ee}{\end{equation}}
\newcommand{\eea}{\end{eqnarray}}
\newcommand{\epr}{\end{proof}}
\newcommand{\ec}{\end{cor}}
\newcommand{\eco}{\end{cond}}
\newcommand{\eex}{\end{ex}}

\newcommand{\nn}{ \nonumber \\ }
\newcommand{\pr}{{\em Proof:\ }}
\newcommand{\qed}{\vrule height 5pt width 5pt depth 0pt}

\newcommand{\err}{\Upsilon}
\newcommand{\id}{{\rm id}}
\newcommand{\p}{{\rm p}} 
\newcommand{\ui}{^{(1)}}
\newcommand{\uii}{^{(2)}}
\newcommand{\di}{_{(1)}}
\newcommand{\dii}{{}_{(2)}}
\newcommand{\ot}{\otimes}
\newcommand{\asr}{{A^R}}
\newcommand{\atr}{{^R\! A}}
\newcommand{\asl}{{_L\! A}}
\newcommand{\atl}{{A_L}}
\newcommand{\asrd}{{{\cal A}^*}}
\newcommand{\atrd}{{^*\!\! {\cal A}}}
\newcommand{\atld}{{{\cal A}_*}}
\newcommand{\asld}{{_*{\cal A}}}
\newcommand{\ci}{\circ}
\newcommand{\ila}{{\cal I}^L({\cal A})}

\newcommand{\fsr}{{\ell_R}}
\newcommand{\ftr}{{_R\ell}}
\newcommand{\fsl}{{_L\ell}}
\newcommand{\ftl}{{\ell_L}}
\newcommand{\lu}{\leftharpoonup}
\newcommand{\ld}{\leftharpoondown}
\newcommand{\ru}{\rightharpoonup}
\newcommand{\rd}{\rightharpoondown}
\newcommand{\lsr}{\lambda\!^*}
\newcommand{\ltr}{^*\!\! \lambda}
\newcommand{\phisl}{_* \!\phi}
\newcommand{\phisr}{\phi\!^*}
\newcommand{\phitr}{^*\!\! \phi}
\newcommand{\phitl}{\phi_*}
\newcommand{\End}{{\rm End}} 
\newcommand{\conv}{\ast}
\def\under{\mbox{\rm\_}\,} 
\newcommand{\F}{\mathcal{F}}
\def\op{^{\rm op}}
\newcommand{\C}{\mathnormal{C}}
\newcommand{\rarr}[1]{\stackrel{#1}{\longrightarrow}} 
\newcommand{\larr}[1]{\stackrel{#1}{\longleftarrow}} 
\newcommand{\M}{\mathcal{M}}
\def\iso{\rarr{\sim}}
\def\Hom{{\rm Hom}} 

\newcommand{\text}{\mbox}
\newcommand{\V}{\mathcal{V}}
\newcommand{\Mod}{\mathbf{Mod}}

\newcommand{\Cat}{\mathsf{Cat}}
\newcommand{\bra}{\langle}
\newcommand{\ket}{\rangle}

\newcommand{\setc}[1]{\setcounter{equation}{#1}}

\newcommand{\lb}{\label}

\begin{document}

 
\large
\title{\bf Hopf algebroids with bijective antipodes:\\ 
           axioms, integrals and duals }
 
\author{\sc Gabriella B\"ohm $^{1,2}$ and
            Korn\'el Szlach\'anyi $^{1,3}$ }

\date{}
 
\maketitle
\normalsize 

\footnotetext[1]{
Research Institute for Particle and Nuclear Physics, Budapest,
H-1525 Budapest 114, P.O.B. 49, Hungary}
\footnotetext[2]{
E-mail: BGABR@rmki.kfki.hu\\
  Supported by the Hungarian Scientific Research Fund, OTKA --
  T 020 285, FKFP -- 0043/2001 and the Bolyai J\'anos Fellowship}
 

\footnotetext[3]{
E-mail: SZLACH@rmki.kfki.hu}

\vskip 2truecm

\begin{abstract}  
Motivated by the study of depth 2 Frobenius extensions  
we introduce a new notion of Hopf algebroid. 
It is a 2-sided bialgebroid with a bijective antipode 
which connects the two, left and right handed, structures.
While all the interesting examples of the Hopf
algebroid of J.H. Lu turn out to be Hopf algebroids in the sense of this
paper, there exist simple examples showing that our definition is
not a special case of Lu's. Our Hopf algebroids, however, belong to the 
class of $\times_L$-Hopf algebras proposed by P. Schauenburg. 
After discussing the axioms and some examples we study the theory of
non-degenerate integrals in order to obtain duals of Hopf algebroids.      
\end{abstract}
\bigskip
\bigskip
\bigskip
\bigskip

\normalsize 


\section{Introduction}

There is a consensus in the literature that bialgebroids, invented by
Takeuchi \cite{T} as $\times_R$-bialgebras, are the proper
generalizations of bialgebras to non-commutative base rings
\cite{Lu,Xu,Schauenburg: Bial,BreMi,Sz,Sz2}. The situation
of Hopf algebroids, i.e., bialgebroids with some sort of antipode, is less
understood. The antipode proposed by J. H. Lu \cite{Lu} is burdened by the 
need of a section for the canonical epimorphism $A\ot A\to A\ot_R A$ the
precise role of which remained unclear. The $\times_R$-Hopf algebras
proposed by P. Schauenburg in \cite{Sch} have a clearcut categorical meaning. 
They are the bialgebroids $A$ over $R$ such that the forgetful functor
$_A\M\to\,_R\M_R$ is not only strict monoidal, which is the condition for $A$
to be a bialgebroid over $R$, but preserves the closed structure as well. In
this very general quantum groupoid, however, antipode, as a map $A\to A$, does
not exist. 

Our proposal of an antipode, announced in \cite{B}, is based on the following
simple observation. The antipode of a Hopf algebra $H$ is a bialgebra map
$S\colon H\to H\op_{\rm cop}$. The opposite of a bialgebroid $A$, however, is
not a bialgebroid in the same sense. In the terminology of \cite{KSz} there
are left bialgebroids and right bialgebroids; corresponding to whether $_A\M$
or $\M_A$ is given a monoidal structure. This suggests that the existence of
antipode on a bialgebroid should be accompanied with a two-sided bialgebroid
structure and the antipode should swap the left and right handed structures.
More explicit guesses for what to take as a definition of the antipode can
be obtained by studying depth 2 Frobenius extensions.  In \cite{KSz}
it has been
shown that for a depth 2 ring extension $N\subset M$ the endomorphism ring
$A=\End\,_NM_N$ has a canonical left bialgebroid structure over the
centralizer $R={\rm C}_M(N)$. If $N\subset M$ is also Frobenius and a
Frobenius homomorphism $\psi\colon M\to N$ is given then $A$ has a right
bialgebroid structure, too. There is a candidate for the antipode $S\colon
A\to A$ as transposition w.r.t. the bilinear form $m,m'\in M\mapsto
\psi(mm')$, see (\ref{eq: S as transpose}). In fact, this definition of $S$
does not require the depth 2 property, so in this example antipode exists
prior to comultiplication.  Assuming the extension $M/N$ is either H-separable
or Hopf-Galois L. Kadison has shown \cite{K1,K2} that the bialgebroid $A$ or
its dual $B=(M\ot_N M)^N$ has an antipode in the sense of \cite{Lu}. 

In a recent paper \cite{DS} B. Day and R. Street give a new
characterization of bialgebroids in the framework of symmetric
monoidal autonomous bicategories. They also introduce a new notion
called {\em Hopf bialgebroid}. It is more restrictive than
Schauenburg's $\times_L$-Hopf algebra since it requires star autonomy
\cite{Barr} rather than closedness. 
We will show in Subsection \ref{DS} that the categorical
definition \cite{DS} of the Hopf bialgebroid is equivalent
to our purely algebraic Definition \ref{def} of {\em Hopf algebroid},
apart from the tiny difference that we allow $S^2$ to be nontrivial on the
base ring. This freedom can be adjusted to the Nakayama automorphism of $\psi$
in case of the Frobenius depth 2 extensions. It is a
new feature of our Hopf algebroids, compared to (weak) Hopf algebras, that
the antipode is not unique. The various antipodes on a given bialgebroid were
shown to be in one-to-one correspondence with the generalized characters
(called twists) in \cite{B}.

It is encouraging that one can find `quantum groupoids' in the
literature that satisfy our axioms. Such are the weak Hopf algebras with
bijective antipode, the
examples of Lu-Hopf algebroids in \cite{BreMi} and the extended Hopf
algebras in \cite{KhalRang}. In particular the Connes-Moscivici
algebra \cite{ConnMo} is a Hopf algebroid it this sense. We also
present an example which is a 
Hopf algebroid in the sense of this paper but does not satisfy the
axioms of \cite{Lu}. This proves that the two notions of Hopf
algebroid -- the one in the sense of this paper and the one in the
sense of \cite{Lu} -- are not equivalent. Until now we could neither
prove nor exclude by examples the possibility that the latter was a
special case of the former. 

The left and right {\em integrals} in Hopf algebra theory are introduced as
the invariants of the left and right regular module, respectively.
In this analogy one can define left integrals in a left bialgebroid
and right integrals in a right bialgebroid. Since a Hopf algebroid has
both left and right bialgebroid structures both left and right
integrals can be defined. 

The properties of the integrals in a (weak) Hopf algebra over a
commutative ring $k$ carry information about its algebraic
structure. For example the Maschke's theorem \cite{LaSw,BNSz} states
that it is a semi-simple algebra if and only if it has a normalized
integral. The Larson-Sweedler theorem \cite{LaSw,V} implies that it is
finite dimensional over $k$ if and only if it has a non-degenerate
left (hence also a right) integral. In this case the $k$-dual also has
a (weak) Hopf algebra structure. 

Therefore, in Section \ref{int}, we analyze the
consequences of the existence of a non-degenerate integral in a Hopf
algebroid. We show that if  there exists a non-degenerate integral in
 a Hopf algebroid ${\cal A}$ over the base $L$ then the ring extension
$L\to A$ is a Frobenius extension hence also finitely generated
projective. We do {\em not} investigate, however, the opposite
implication i.e. we do not study the question under what conditions on
the Hopf algebroid the existence of a non-degenerate integral
follows. 

If some of the $L$-module structures of a Hopf algebroid is
finitely generated projective then the corresponding dual can be
equipped with a bialgebroid structure \cite{KSz}. 
There is no obvious way, however, how to equip it with an antipode in
general.  
We show that in the case of Hopf algebroids possessing a
non-degenerate integral the  dual
bialgebroids are all (anti-) isomorphic, and they combine into a Hopf
algebroid -- depending on the choice of the non-degenerate integral.
Therefore we do not associate a dual Hopf algebroid to a given Hopf
algebroid rather a dual isomorphism class to an isomorphism class of
Hopf algebroids. The well known $k$-dual of a finite (weak) Hopf
algebra $H$ over the commutative ring $k$ turns out to be the
unique (distinguished) (weak) Hopf algebra in the dual isomorphism
class of the isomorphism class of the Hopf algebroid $H$. 

The paper is organized as follows: In Section \ref{preli} we
introduce some technical conventions about bialgebroids that are used in this
paper. Our motivating example, the Hopf algebroid
corresponding to a depth 2 Frobenius extension of rings is discussed
in Section \ref{ringex}. In Section \ref{defsec} we give
some equivalent definitions of Hopf algebroids. We prove that our
definition is equivalent to the one in \cite{DS} hence gives a special
case of the one in \cite{Sch}. In the final subsection of Section
\ref{defsec} we present a collection of examples. In the Section
\ref{int} we propose a theory of non-degenerate integrals as a tool for the
definition of the dual Hopf algebroid.

\section{Preliminaries on  bialgebroids}
\lb{preli}
\setc{0}
In this technical section we summarize our notations and the basic definitions
of bialgebroids that will be used later on. For more about
bialgebroids we refer to the literature \cite{T,Sch,BreMi,KSz,Sch1,Sz2,Sz: EM}.

\bd \lb{lbgd}
A {\em left bialgebroid} (or Takeuchi $\times_L$-bialgebra) ${\cal
A}_L$ consists of the data $(A,L,s_L,t_L,$\break$\gamma_L,\pi_L)$. The
$A$ and 
$L$ are associative unital rings, the total and base rings,
respectively. The $s_L:L\to A$ and $t_L:L^{op}\to 
A$ are ring homomorphisms such that the images of $L$ in $A$ commute
making $A$ an $L$-$L$ bimodule via
\be  l\cdot a\cdot l^{\prime}\colon = s_L(l) t_L(l^{\prime}) a .
\lb{elbim} \ee
The bimodule (\ref{elbim}) is denoted by ${_L A_L}$. The triple
$({_L A_L},\gamma_L,\pi_L)$ is a comonoid in ${_L {\cal M}_L}$, the
category of $L$-$L$-bimodules. Introducing
the Sweedler's notation $\gamma_L(a)\equiv a\di \ot a\dii\in \atl \ot
\asl$ the identities
\bea  a\di t_L(l) \ot a\dii &=& a\di \ot a\dii s_L(l) \lb{cros}\\
      \gamma_L(1_A)&=& 1_A\ot 1_A \\
      \gamma_L(ab)&=&\gamma_L(a) \gamma_L(b) \lb{gmp} \\
      \pi_L(1_A) &=& 1_L \\
      \pi_L\left(a s_L\ci \pi_L(b)\right)=&\pi_L&(ab)=
       \pi_L\left(a t_L\ci \pi_L(b)\right)
\eea
are  required for all $l\in L$ and $a,b\in A$. The requirement 
(\ref{gmp}) makes sense in the view of (\ref{cros}).
\ed
The $L$ actions of the bimodule ${_L\!A_L}$ in (\ref{elbim}) are given
by left multiplication. Using right multiplication there exists
another $L$-$L$ bimodule structure on the total ring $A$ of a left
bialgebroid ${\cal A}_L$:
\be l\cdot a\cdot l^{\prime} \colon = a t_L(l)
s_L(l^{\prime}). \lb{elbimr}\ee 
This $L$-$L$ bimodule is called ${^L\!A^L}$. This way $A$ carries four
commuting actions of $L$.

If ${\cal A}_L=(A,L,s_L,t_L,\gamma_L,\pi_L)$ is a left bialgebroid
then so is its co-opposite: ${\cal
A}_{L\ cop}=$\break $(A,L^{op},t_L,s_L,\gamma_L^{op},\pi_L)$. The
opposite ${\cal 
A}_L^{op}=(A^{op},L,t_L,s_L,\gamma_L,\pi_L)$ has a different
structure that was introduced under the name {\em right bialgebroid}
in \cite{KSz}: 
\bd A {\em right bialgebroid} ${\cal A}_R$ consists of the data
$(A,R,s_R,t_R,\gamma_R,\pi_R)$. The $A$ and 
$R$ are associative unital rings, the total and base rings,
respectively. The $s_R:R\to A$ and $t_R:R^{op}\to 
A$ are ring homomorphisms such that the images of $R$ in $A$ commute
making $A$ an $R$-$R$ bimodule:
\be  r\cdot a\cdot r^{\prime}\colon = a s_R(r^{\prime}) t_R(r^).
\lb{erbim}\ee
The bimodule (\ref{erbim}) is denoted by ${^R\! A^R}$. The triple
$({^R\! A^R},\gamma_R,\pi_R)$ is a comonoid in ${_R {\cal
M}_R}$. Introducing 
the Sweedler's notation $\gamma_R(a)\equiv a\ui \ot a\uii\in \asr \ot
\atr$ the identities
\bea  s_R(r)a\ui  \ot a\uii &=& a\ui \ot  t_R(r) a\uii \nn
      \gamma_R(1_A)&=& 1_A\ot 1_A \nn
      \gamma_R(ab)&=&\gamma_R(a) \gamma_R(b) \nn
      \pi_R(1_A) &=& 1_R \nn
      \pi_R\left(s_R\ci \pi_R(a) b \right)=&\pi_R&(ab)=
       \pi_R\left(t_R\ci \pi_R(a) b \right)
\nonumber\eea
are  required for all $r\in R$ and $a,b\in A$.
\ed
For the right bialgebroid ${\cal A}_R$ we introduce the
$R$-$R$ bimodule ${_R\! A_R}$ via
\be  r\cdot a\cdot r^{\prime}\colon = s_R(r) t_R(r^{\prime})
a. \lb{erbiml}\ee 
This way $A$ caries four commuting actions of $R$.

Left (right) bialgebroids can be characterized by the property that the
forgetful functor ${_A{\cal M}}\to {_L{\cal M}_L}$ ( ${{\cal
M}_A}\to {_R{\cal M}_R}$ ) is strong monoidal \cite{Schauenburg: Bial,Sz}.

It is natural to consider the homomorphisms of bialgebroids to be ring
homomorphisms preserving the comonoid structure. We do not want to make
difference however between bialgebroids over isomorphic base
rings. This leads to the 
\bd \cite{Sz2} A left bialgebroid homomorphism ${\cal A}_L \to {\cal
A}^{\prime}_L$ is a pair of ring homomorphisms $(\Phi:A\to A^{\prime},
\phi: L\to L^{\prime})$ such that 
\bea s^{\prime}_L\ci \phi &=& \Phi\ci s_L\nn
     t^{\prime}_L\ci \phi &=& \Phi\ci t_L\nn
     \pi^{\prime}_L\ci \Phi &=& \phi\ci \pi_L\nn
     \gamma_L^{\prime}\ci \Phi&=& (\Phi\ot \Phi)\ci \gamma_L.
\nonumber\eea
The last condition makes sense since by the first two
conditions $\Phi\ot \Phi$ is a well defined map $\atl \ot \asl\to
{A_{L^{\prime}}}\ot {_{L^{\prime}}}A$.

The pair $(\Phi,\phi)$ is an isomorphism of left bialgebroids if it is
a bialgebroid homomorphism such that both $\Phi$ and $\phi$ are 
bijective. 

A right bialgebroid homomorphism (isomorphism) ${\cal A}_R\to {\cal
A}^{\prime}_R$ is a left bialgebroid homomorphism (isomorphism) $({\cal
A}_R)^{op}\to ({\cal A}^{\prime}_R)^{op}$.
\ed
Let ${\cal A}_L$ be a left bialgebroid. The equation (\ref{elbim})
describes two $L$-modules $\atl$ and $\asl$. Their $L$-duals
are the additive groups of $L$-module maps:
\[  \atld\colon = \{ \phitl: \atl \to {L_ L} \} \quad {\rm and} \quad
     \asld\colon = \{ \phisl: \asl \to {_L L} \} \]
where ${_L L}$ stands for the left regular and $L_L$ for the right
regular $L$-module. Both $\atld$ and $\asld$ carry left $A$ module
structures via the transpose of the right regular action of $A$. For
$\phitl\in \atld, \phisl\in\asld$ and $a,b\in A$ we have:
\[   \left(a\ru \phitl\right)(b) =\phitl(ba)\quad {\rm and}\quad
     \left(a\rd \phisl\right)(b) =\phisl(ba). \]
Similarly, in the case of a right bialgebroid ${\cal A}_R$ --  denoting
the left and right regular $R$-modules by $^R R$ 
and $R^R$, respectively, -- the two $R$-dual additive groups
\[  \asrd\colon = \{ \phisr: \asr \to {R^R} \} \quad {\rm and}\quad
    \atrd\colon = \{ \phitr: \atr \to {^R R} \} \]
carry right $A$-module structures:
\[ \left(\phisr \lu a \right)(b) =\phisr(ab) \quad {\rm and}\quad
   \left(\phitr \ld a \right)(b) =\phitr(ab). \]
The comonoid structures can be transposed to give monoid (i.e. ring)
structures to the duals. In the case of a left bialgebroid ${\cal A}_L$ 
\be  \left(\phitl {\psi_*} \right)(a)=
{\psi_*}\left( s_L \ci \phitl (a\di) a\dii \right)\quad {\rm and}\quad
\left(\phisl {_* \psi} \right)(a)=
{_* \psi}\left( t_L \ci \phisl (a\dii) a\di \right) \ee
for $\phisl, {_*\psi}\in \asld$, $\phitl,{\psi_*}\in \atld$ and $a\in
A$. 

Similarly, in the case of a right bialgebroid ${\cal A}_R$ 
\be \left(\phisr {\psi^*} \right)(a)=
\phisr\left( a\uii t_R \ci {\psi^*} (a\ui)\right)\quad {\rm and}\quad
 \left(\phitr {^*\! \psi} \right)(a)=
\phitr \left(a\ui  s_R \ci {^* \!\psi} (a\uii)\right)\ee
for $\phisr, {\psi^*}\in \asrd$, $\phitr,{^*\psi}\in \atrd$ and $a\in
A$.

In the case of a left bialgebroid ${\cal A}_L$  also the ring $A$ has 
right $\atld$- and   right $\asld$- module structures:
\be a\lu \phitl = s_L\ci \phitl (a\di)a\dii \quad {\rm and}\quad
     a\ld \phisl = t_L\ci \phisl (a\dii)a\di \ee
for $\phitl\in \atld$, $\phisl\in \asld$ and $a\in A$.

Similarly, in the case of a right bialgebroid ${\cal A}_R$  the ring
$A$ has left $\asrd$- and left $\atrd$ structures:
\be \phisr\ru a = a\uii t_R\ci \phisr(a\ui )\quad {\rm and}\quad
     \phitr\rd a = a\ui  s_R\ci \phitr(a\uii) \lb{srac}\ee
for $\phisr\in \asrd$,  $\phitr\in \atrd$ and $a\in A$.

In the case when the $L$ ($R$) module structure on $A$ is finitely
generated projective then the corresponding dual has also a
bialgebroid structure: If  ${\cal A}_L$ is a left bialgebroid such
that the $L$-module $\atl$ is finitely generated projective then
$\atld$ is a right bialgebroid over the base ${L_*}\equiv L$ as follows:
\[ \left({s_{*R}}(l) \right)(a) =\pi_L(a s_L(l)), \quad
     \left({t_{*R}}(l) \right)(a) = l \pi_L(a) ,\quad
     \gamma_{*R} (\phitl) = b_i \ru \phitl \ot {\beta^i_*},\quad
     \pi_{*R}(\phitl) = \phitl(1_A) \]
where $\{b_i\}$ is an $L$-basis in $\atl$ and $\{\beta^i_*\}$ is the dual
basis in $\atld$ . 

Similarly, if  ${\cal A}_L$ is a left bialgebroid such
that the $L$-module $\asl$ is finitely generated projective then
$\asld$ is a right bialgebroid over the base ${_*\! L}\equiv L$ as
follows: 
\[ \left({_*s_{R}}(l) \right)(a) = \pi_L(a) l,\quad
     \left({_*t_{R}}(l) \right)(a) = \pi_L(a t_L(l)), \quad
     {_*\gamma_{R}} (\phisl) = {_*\beta^i}\ot b_i \rd \phisl, \quad
     {_*\pi_{R}}(\phisl)=  \phisl(1_A)\]
where $\{b_i\}$ is an $L$-basis in $\asl$ and $\{_*\beta^i\}$ is the dual
basis in $\asld$ . 

If  ${\cal A}_R$ is a right bialgebroid such
that the $R$-module $\asr$ is finitely generated projective then
$\asrd$ is a left bialgebroid over the base $R^*\equiv R$ as follows:
\[ \left({s^*_{L}}(r) \right)(a) = r\pi_R(a),\quad
     \left({t^*_{L}}(r) \right)(a) =  \pi_R(s_R(r)a),\quad
     {\gamma^*_{R}} (\phisr) = \phisr\lu b_i \ot {\beta^{*i}},\quad
     {\pi^*_{R}}(\phisr) = \phisr(1_A)\]
where $\{b_i\}$ is an $R$-basis in $\asr$ and $\{\beta^{*i}\}$ is the dual
basis in $\asrd$ . 

If  ${\cal A}_R$ is a right bialgebroid such
that the $R$-module $\atr$ is finitely generated projective then
$\atrd$ is a left bialgebroid over the base ${^*\!R}\equiv R$ as follows:
\[ \left({^*s_{L}}(r) \right)(a) = \pi_R(t_R(r)a),\quad
     \left({^*t_{L}}(r) \right)(a) = \pi_R(a) r ,\quad
     {^*\gamma_{R}} (\phitr) = {^*\beta^i}\ot \phitr\ld b_i ,\quad
     {^*\pi_{R}}(\phitr) = \phitr(1_A)\]
where $\{b_i\}$ is an $R$-basis in $\atr$ and $\{^*\beta^{i}\}$ is the dual
basis in $\atrd$ . 


\section{The motivating example: D2 Frobenius extensions}
\lb{ringex}
\setc{0}

\subsection{The forefather of antipodes} 

In this subsection $N\to M$ denotes a Frobenius extension of rings. 
This means the existence of $N$-$N$ bimodule
maps $\psi\colon M\to N$ possessing quasibases. An element $\sum_i
u_i\ot v_i\in 
M\stackrel{\ot}{_{_N}} M$ is called the quasibasis of $\psi$ \cite{W} if 
\be
\sum_i \psi(mu_i)\cdot v_i\ =\ m\ =\ \sum_i u_i\cdot\psi(v_i m)\,,\qquad m\in
M\,.
\lb{quasib}\ee
As we shall see, already in this general
situation there exist antiautomorphisms $S$ on the ring
$A:=\End\,_NM_N$, one for each Frobenius homomorphism $\psi$. The $S$ will
become an antipode if the extension $N\subset M$ is also of depth 2, so $A$
also has coproduct(s).

In addition to $A$, important role is played by the center of the bimodule
$M\stackrel{\ot}{_{_N}} M$ 
\[
B:=(M\stackrel{\ot}{_{_N}} M)^N\equiv\{X\in M\stackrel{\ot}
{_{_N}}M\,|\,n\cdot X=X\cdot n\ \forall n\in N\,\}
\]
which is a ring with multiplication $(b^1\ot b^2)({b'}^1\ot {b'}^2)={b'}^1b^1\ot
b^2{b'}^2$ and unit $1_B=1_M\ot 1_M$. Note that the ring  structures of neither
$A$ nor $B$ depend on the Frobenius structure. But if there is a Frobenius
homomorphism $\psi$ then Fourier transformation makes $A$ and $B$ isomorphic
as additive groups.

Fixing a Frobenius homomorphism $\psi$ with quasibasis $\sum_i u_i\ot v_i$ we
can introduce convolution products on both $A$ and $B$ as follows. From now on
we omit the summation symbol for summing over the quasibasis. 
\bea
\alpha,\beta\in
A&\mapsto&\alpha\conv\beta:=\alpha(u_i)\beta(v_i\under)\ 
\in A\nn
a,b\in B&\mapsto& a\conv b:=a^1\psi(a^2b^1)\ot b^2\ \in B\,.
\nonumber\eea
The convolution product lends $A$ and $B$ new ring structures. The unit
of $A$ is $\psi$ and the unit of $B$ is $u_i\ot v_i$. The Fourier
transformation is to relate these new algebra structures to the old ones.
There are two natural candidates for a Fourier transformation,
\bea
{\F} \colon A\to B\,,\quad&{\F}(\alpha):=u_i\ot\alpha(v_i)\nn
\dot{\F}\colon A\to B\,,\quad&\dot{\F}(\alpha):=\alpha(u_i)\ot v_i\nn
{\F}^{-1}\colon B\to A\,,\quad&{\F}^{-1}(b)=\psi(\under b^1)b^2\nn
\dot{\F}^{-1}\colon B\to A\,,\quad&\dot{\F}^{-1}(b)=b^1\psi(b^2\under)
\nonumber\eea
They relate the convolution and ordinary products or their opposites
as follows. 
\bea
{\F}(\alpha\conv\beta)&=&{\F}(\alpha){\F}(\beta)\nn
\dot{\F}(\alpha\conv\beta)&=&\dot{\F}(\beta)\dot{\F}(\alpha)\nn
{\F}(\alpha\beta)&=&{\F}(\beta)\conv{\F}(\alpha)\nn
\dot{\F}(\alpha\beta)&=&\dot{\F}(\alpha)\conv\dot{\F}(\beta)
\nonumber\eea
The difference between ${\F}$ and $\dot{\F}$ is therefore an antiautomorphism on
both $A$ and $B$. This leads to the "antipodes"
\bea
S_A\colon A\to A\op\,,\quad&S_A:=\dot{\F}^{-1}\circ{\F}\,,\quad &
S_A(\alpha)=u_i\psi(\alpha(v_i)\under) \lb{eq: S_A}\\
S_B\colon B\to B\op\,,\quad&S_B:=\dot{\F}\circ{\F}^{-1}\,,\quad &
S_B(b)=\psi(u_ib^1)b^2\ot v_i
\eea
with inverses
\bea
S_A^{-1}(\alpha)&=&\psi(\under\alpha(u_i))v_i \lb{eq: S inv}\\
S_B^{-1}(b)&=&u_i\ot b^1\psi(b^2v_i)
\eea
Notice that $S_A$ is just transposition w.r.t. the bi-$N$-linear form
$(m,m')=\psi(mm')$ since
\be \lb{eq: S as transpose}
\psi(mS_A(\alpha)(m'))=\psi(\alpha(m)m')\,.
\ee

These antipodes behave well also relative to the 
bimodule structures over the centralizer $\C_M(N):=\{c\in M|cn=nc,\ \forall
n\in N\}$. Let us consider $S_A$. The centralizer is embedded into $A$ twice:
via left multiplications and right multiplications, 
$$
L\rarr{\lambda}\ A\ \larr{\rho}R
$$
where $L$ stands for $\C_M(N)$ and $R$ for $\C_M(N)\op$. 
Clearly, $\lambda(L)\subset\C_M(\rho(R))$. Introducing the
Nakayama automorphism
$$
\nu\colon \C_M(N)\to\C_M(N)\,,\qquad \nu(c):=\psi(u_ic)v_i
$$ 
of $\psi$ and using its basic identities
\bea
\psi(mc)&=&\psi(\nu(c)m)\,,\qquad m\in M,\ c\in\C_M(N)\nn
u_ic\ot v_i&=&u_i\ot\nu(c) v_i\,,\qquad c\in\C_M(N)
\nonumber\eea
we obtain
\bea
S_A\circ\lambda&=&\rho\circ\nu^{-1}\nn
S_A\circ\rho&=&\lambda\nonumber
\eea
and therefore
\be \lb{eq: .S.}
S_A(\lambda(l)\rho(r)\alpha)=S_A(\alpha)\lambda(r)\rho(\nu^{-1}(l))\,,\qquad
l\in L,\ r\in R,\ \alpha\in A\,.
\ee
In order to interpret the latter relation as the statement that $S_A$
is a bimodule map we define $L$-$L$ and $R$-$R$ bimodule
structures on $A$ by
\bea
l_1\cdot\alpha\cdot l_2&:=&s_L(l_1)t_L(l_2)\alpha \lb{eq: lAl}\\
r_1\cdot \alpha\cdot r_2&:=&\alpha t_R(r_1)s_R(r_2) \lb{eq: rAr}
\eea
where we introduced the ring homomorphisms
\be
\begin{array}{ll}
s_L:=L\rarr{\lambda} A\qquad &s_R:=R\rarr{\rho}A \lb{eq: s}\\
t_L:=L\op\rarr{\id}R\rarr{\rho}A\qquad &t_R:=R\op\rarr{\nu}L\rarr{\lambda}A
\lb{eq: t}
\end{array}
\ee
Also using the notation $\theta$ for the inverse of the Nakayama automorphism
when considered as a map
\[
\theta\colon L\rarr{\nu^{-1}}R\op
\]
equation (\ref{eq: .S.}) can be read as 
\be \lb{eq: rrSll}
S_A(l_1\cdot \alpha\cdot l_2)\ =\
\theta(l_2)\cdot S_A(\alpha)\cdot\theta(l_1)\,.
\ee

\br
The apparent asymmetry between $t_L$ and $t_R$ in (\ref{eq: t}) disappears if
one repeats the above construction for the more general situation of a
Frobenius $N$-$M$ bimodule $X$ instead of the
$_NM_M$ arising from a Frobenius extension of rings. As a matter
of fact, denoting by $\bar X$ the (two-sided) dual of $X$ and
setting $A=\End\, X\stackrel{\ot}{_{M}}\bar X$, $L=\End X$ and $R=\End
\bar X$ we find 
the obvious ring homomorphisms
$$
s_L(l)=l\ot \bar X\,,\qquad s_R(r)=X\ot r 
$$
but there is no distinguished map $L\to R\op$ like the identity is in the case
of $X=\,_NM_M$. Instead we have two distinguished maps given by the left and
right dual functors (transpositions). It is easy to check that in case of
$X=\,_NM_M$ they are the identity and the $\nu^{-1}$, respectively, as we used
in (\ref{eq: t}). 
\er

In order to restore the symmetry let us introduce the counterpart of $\theta$
which is the identity as a homomorphism $\iota\colon R\rarr{\id}L\op$. Then,
in addition to (\ref{eq: rrSll}) the antipode satisfies also
\be \lb{eq: llSrr}
S_A(r_1\cdot\alpha\cdot r_2)\ =\ 
\iota(r_2)\cdot S_A(\alpha)\cdot \iota(r_1)\,.
\ee 
The most important consequence of (\ref{eq: rrSll}) and (\ref{eq: llSrr}) is
the existence of a tensor square of $S_A$. In the case of Hopf algebras one
often uses expressions like $(S\ot S)\circ\Sigma$, where $\Sigma$ is
the symmetry $A\ot A\to A\ot A$, $x\ot y\mapsto y\ot x$ in the category of
$k$-modules. Now we  have bimodule categories ${_L\M_L}$ and ${_R\M_R}$
without braiding so $\Sigma$ does 
not exist and neither do $S_A\ot S_A$ or $S_A^{-1}\ot S_A^{-1}$ because $S_A$ is
not a bimodule map. Instead we have the \textit{twisted bimodule properties}
(\ref{eq: rrSll})  and (\ref{eq: llSrr}) which guarantee the existence of the
`composite of' $S_A\ot S_A$ and $\Sigma$ although individually they don't exist.
More precisely, there exist twisted bimodule maps
\bea
S_{A\stackrel{\ot}{_{_L}} A}\colon A\stackrel{\ot}{_{_L}} A&\to& A\stackrel{\ot}
{_{_R}} A\nn 
\alpha\stackrel{\ot}{_{_L}}\beta&\mapsto& S_A(\beta)\stackrel{\ot}
{_{_R}} S_A(\alpha)\nn
S_{A\stackrel{\ot}{_{_R}} A}\colon A\stackrel{\ot}{_{_R}} A&\to& A\stackrel{\ot}
{_{_L}} A\nn
\alpha\stackrel{\ot}{_{_R}}\beta&\mapsto&  S_A(\beta)\stackrel{\ot}{_{_L}} 
S_A(\alpha)\,.
\nonumber\eea

For later convenience let us record some useful formulas following directly
from (\ref{eq: S_A}) and (\ref{eq: t}).
\be\begin{array}{ll}
S_A\circ s_L=t_L\circ\nu^{-1}&\qquad
S_A\circ t_L=s_L\nn
S_A\circ s_R=t_R\circ\nu^{-1}&\qquad
S_A\circ t_R=s_R\nn
t_L\circ\iota=s_R&\qquad
t_R\circ\theta=s_L\lb{eq:Sbimm}
\end{array}\ee
Notice also that $s_L(L)$ and $t_R(R)$ are the same subrings of $A$ and
similarly $t_L(L)=s_R(R)$.

\subsection{Two-sided bialgebroids}

Recall from \cite{KSz} that for depth 2 extensions $N\to M$ the $A$
has a canonical left bialgebroid structure over $L$ in which the
coproduct $\gamma_L\colon A\to A\stackrel{\ot}{_{_L}} 
A$ is an $L$-$L$ bimodule map with respect to the bimodule structure
(\ref{eq: lAl}). If $N\to M$ is also Frobenius then there is another right
bialgebroid structure on $A$, canonically associated to a choice of $\psi$, in
which $R$ is the base and $A$ is an $R$-$R$ bimodule via (\ref{eq: rAr}).
Moreover these two structures are related by the antipode. This two-sided
structure is our motivating example of a Hopf algebroid.

We start with a technical lemma on the left and right quasibases. $A$ and $B$
denotes the rings as before.  
\bl \lb{lem: qbases}
Let $N\to M$ be a Frobenius extension and $\psi$, $u_i\ot v_i$ be a fixed
Frobenius structure. Let $n$ be a positive integer and 
$\beta_i$, $\gamma_i\in A$ and $b_i$, $c_i\in B$, for $i=1,\dots, n$. Assume
they are related via $b_i={\F}(\gamma_i)$ and $c_i=\dot{\F}(\beta_i)$. Then the
following conditions are equivalent (summation symbols over $i$ suppressed):
\[\begin{array}{rll}
i)&b^1_i\stackrel{\ot}{_{_N}} b^2_i\beta_i(m)\ =\ m\stackrel{\ot}{_{_N}} 
1_M\,,&\qquad m\in M\nn
ii)&\gamma_i(m)c_i^1\stackrel{\ot}{_{_N}} c_i^2\ =\ 1_M\stackrel{\ot}{_{_N}} 
m\,,&\qquad m\in M\nn
iii)&\gamma_i(m)\beta_i(m')\ =\ \psi(mm')\,,&\qquad m,m'\in M\nn
iv)&b_i^1\stackrel{\ot}{_{_N}} b_i^2c_i^1\stackrel{\ot}{_{_N}} c_i^2\ 
=\ u_k\stackrel{\ot}{_{_N}} 1_M\stackrel{\ot}{_{_N}} v_k& \nonumber
\end{array}\]
\el
If such elements exist the extension is called $D2$, i.e., of depth 2.
The first two conditions are meaningful also in the non-Frobenius case and
therefore $\{b_i,\beta_i\}$ was called in \cite{KSz} a left D2 quasibasis and
$\{c_i,\gamma_i\}$ a right D2 quasibasis. The equivalence of conditions i)
and ii) was shown in \cite{KSz}, Proposition 6.4. The rest of the proof is left to
the reader. 

For D2 extensions the map $\alpha\ot\beta\mapsto\{m\ot m'\mapsto
\alpha(m)\beta(m')\}$ is an isomorphism
\[
A\stackrel{\ot}{_{_L}} A\iso\Hom_{N-N}(M\stackrel{\ot}{_{_N}} M,M)\,,
\]
see \cite{KSz}, Proposition 3.11. Then the coproduct $\gamma_L:A\to
A\stackrel{\ot}{_{_L}} A$ is the 
unique map $\alpha\mapsto\alpha\di\ot\alpha\di$ which satisfies
\[
\alpha(mm')\ =\ \alpha\di(m)\alpha\di(m')\,.
\]
For D2 Frobenius extensions we can dualize this construction. We have
the isomorphism
\bea
A\stackrel{\ot}{_{_R}} A&\iso&\Hom_{N-N}(M,M\stackrel{\ot}{_{_N}} M)\nn
\alpha\stackrel{\ot}{_{_R}}\beta&\mapsto&\alpha(\under u_i)\stackrel{\ot}
{_{_N}}\beta(v_i)\,.\nonumber
\eea
Then $\gamma_R\colon A\to A\stackrel{\ot}{_{_R}} A$ is defined as the unique map
$\alpha\mapsto\alpha\ui\ot\alpha\uii$ for which
\be \lb{eq: def cop R}
\alpha(m)u_i\stackrel{\ot}{_{_N}} v_i\ =\ \alpha\ui(mu_i)\stackrel{\ot}
{_{_N}}\alpha\uii(v_i)\,.
\ee
Explicit formulas for both coproducts, as well as their counits, are given in
the Corollary below. But even without these formulas we can find out how the
two coproducts are related by the antipode.
\bt
For a D2 Frobenius extension $N\to M$ of rings the endomorphism ring
$A=\End\,_NM_N$ is a left bialgebroid over $L=\C_M(N)$ and a right bialgebroid
over $R=L\op$ such that the antipode defined in (\ref{eq: S_A}) gives rise to
isomorphisms
\bea \lb{eq: S in Bgd}
  (A, R,s_R, t_R,\gamma_R,\pi_R) &\rarr{(S_A,\iota)}&
\ (A,L,s_L,t_L,\gamma_L,\pi_L)^{op}_{cop}\nn
(A,L,s_L,t_L,\gamma_L,\pi_L)  &\rarr{(S_A,\theta)}&
\ (A, R,s_R, t_R,\gamma_R,\pi_R) ^{op}_{cop}
\eea
of left bialgebroids. That is to say,
\bea
S_A\colon A&\to& A\op\ {is\  a\  ring\  isomorphism}\\
S_A\ci s_R&=&s_L\ci\iota, \ S_A\ci t_R=t_L\ci\iota \quad {and}\quad
S_A\ci s_L = s_R\ci\theta,\ S_A\ci t_L=t_R\ci\theta \\ 
\gamma_L\circ S_A&=&S_{A\stackrel{\ot}{_{_R}} A}\circ \gamma_R \quad {and}\quad 
\gamma_R\circ S_A=S_{A\stackrel{\ot}{_{_L}} A}\circ \gamma_L\lb{eq: copScop}\\
\pi_L\circ S_A&=&\iota\ci \pi_R\quad  {and} \quad  \pi_R\circ
S_A=\theta\circ \pi_L. 
\eea

\et
\pr
The left bialgebroid structure of $A$ has been constructed in
\cite{KSz}, Theorem 4.1. The right bialgebroid structure will follow
automatically after 
establishing the four properties of the antipode. The first two have already
been discussed before. In order to prove (\ref{eq: copScop}) recall the
definition (\ref{eq: def cop R}) of $\gamma_R$. Thus (\ref{eq: copScop}) is
equivalent to 
\bea
S_A^{-1}(S_A(\alpha)\di)(mu_i)\stackrel{\ot}{_{_N}} S_A^{-1}(S_A(\alpha)\di)(v_i)
&=&\alpha(m)u_i\stackrel{\ot}{_{_N}} v_i \lb{eq:scop}\\
S_A(S_A^{-1}(\alpha)\di)(mu_i)\stackrel{\ot}{_{_N}} S_A(S_A^{-1}(\alpha)\di)(v_i)
&=&u_i\stackrel{\ot}{_{_N}} v_i\alpha(m) \lb{eq:sicop}
\eea
for all $m\in M$. Expanding the left hand side of (\ref{eq:scop}) then using 
(\ref{eq: S as transpose}), then (\ref{eq: S inv}), then the
definition of $\gamma_L$, and finally (\ref{eq: S as transpose}) again we
obtain  
\bea
&&S_A^{-1}(S_A(\alpha)\di)(mu_i)\psi(S_A^{-1}(S_A(\alpha)\di)(v_i)u_k)\stackrel{\ot}
{_{_N}} v_k\nn
&&=S_A^{-1}(S_A(\alpha)\di)(mu_i)\psi(v_i S_A(\alpha)\di(u_k))\stackrel{\ot}{_{_N}} 
v_k\nn
&&=S_A^{-1}(S_A(\alpha)\di)(mS_A(\alpha)\di(u_k))\stackrel{\ot}{_{_N}} v_k\nn
&&=\psi(mS_A(\alpha)\di(u_k)S_A(\alpha)\di(u_j))v_j\stackrel{\ot}{_{_N}} v_k\nn
&&=\psi(m S_A(\alpha)(u_ku_j))v_j\stackrel{\ot}{_{_N}} v_k\nn
&&=\psi(\alpha(m)u_ku_j)v_j\stackrel{\ot}{_{_N}} v_k\nn
&&=\alpha(m)u_k\stackrel{\ot}{_{_N}} v_k\nonumber
\eea
and analogously (\ref{eq:sicop}).
Now it is easy to see that both $\pi_L\circ S_A$ and $\theta\circ
\pi_L\circ S_A^{-1}$ are counits for $\gamma_R$.
Therefore both are equal to \textit{the} counit $\pi_R$. This finishes
the proof of the 
isomorphisms (\ref{eq: S in Bgd}).
\qed

\bc
Explicit formulas for the left and right bialgebroid structures can be given
using the quasibases of Lemma \ref{lem: qbases} as follows.
\bea
\gamma_L(\alpha)&=&\gamma_i\stackrel{\ot}{_{_L}} c_i^1\alpha(c_i^2\under)\nn
&=&\alpha(\under b_i^1)b_i^2\stackrel{\ot}{_{_L}} \beta_i\nn
&=&\gamma_i\stackrel{\ot}{_{_L}}\beta_i\conv\alpha\nn
&=&\alpha\conv\gamma_i\stackrel{\ot}{_{_L}}\beta_i\nn
\pi_L(\alpha)&=&\alpha(1_M)\nn
\gamma_R(\alpha)&=&\alpha(\under c_i^1)c_i^2\stackrel{\ot}{_{_R}}
\psi(\under\gamma_i(u_k))v_k\nn
&=&u_k\psi(\beta_i(v_k)\under)\stackrel{\ot}{_{_R}} b_i^1 \alpha(b_i^2\under)\nn
&=&\alpha \conv S_A(\beta_i)\stackrel{\ot}{_{_R}} S_A(\gamma_i)\nn
&=&S_A(\beta_i)\stackrel{\ot}{_{_R}} S_A(\gamma_i)\conv\alpha\nn
\pi_R(\alpha)&=& u_i \psi\ci\alpha(v_i)
\nonumber\eea
\ec


\section{Hopf algebroids}
\setc{0}
\lb{defsec}
\subsection{The Definition}

The total ring of a Hopf algebroid carries eight canonical module
structures over the base ring -- modules of the kind (\ref{elbim}),
(\ref{elbimr}), (\ref{erbim}) and (\ref{erbiml}). In this situation the
standard notation for the tensor product of modules, e.g.
$A\stackrel{\ot}{_{_R}}A$, would be ambiguous. In order to avoid any
misunderstandings we therefore put marks on both modules, as in $\asr\ot \atr$
for example, that indicate the module structures taking part in the tensor
product. Other module structures (commuting with those taking part in the
tensor product) are usually unadorned and should be clear from the context.

For coproduts of left bialgebroids we use the Sweedler's notation in
the form $\gamma_L(a)=a\di \ot a\dii$ and of right bialgebroids
$\gamma_R(a)=a\ui\ot a\uii$.

\bd \lb{def} The Hopf algebroid is a pair $({\cal A}_L,S)$ consisting
of a left 
bialgebroid ${\cal A}_L=(A,L,s_L,$ $t_L,\gamma_L,\pi_L)$
and an 
anti-automorphism $S$ of the total ring $A$ satisfying
\bea i)&& S\ci t_L=s_L \qquad {and}\lb{axi}\\
    ii)&& S^{-1}(a\dii)_{(1^{\prime})}\ot S^{-1}(a\dii)_{(2^{\prime})}
    a\di= S^{-1}(a)\ot 1_A \lb{axsi}\\
       && S(a\di)_{(1^{\prime})}a\dii \ot  S(a\di)_{(2^{\prime})}=
    1_A\ot S(a)\lb{axs}
\eea
as elements of $\atl\ot \asl$, for all $a\in A$.
\ed
The axiom (\ref{axs}) implies that 
\be S(a\di)a\dii=t_L\ci \pi_L\ci S(a)\lb{luap}\ee
for all $a\in A$. 
Introduce the map $\theta_L\colon =
\pi_L\ci S\ci s_L :L\to L$. Owing to (\ref{luap}) it satisfies
\bea t_L\ci \theta_L(l)&=& t_L\ci \pi_L \ci S \ci s_L(l) =S\ci s_L(l)\nn
\theta_L(l) \theta_L(l^{\prime})&=& 
\pi_L\ci S\ci s_L(l)\pi_L\ci S\ci s_L (l^{\prime}) = 
\pi_L\left( t_L\ci \pi_L\ci S\ci s_L (l^{\prime})  S\ci s_L(l) \right)=\nn
&=&\pi_L\left(S\ci s_L(l^{\prime})S\ci s_L(l)\right)= \theta_L(l
l^{\prime}).\lb{theta} 
\eea
In view of (\ref{theta}) 
$S$ is a {\em twisted bimodule map} ${^R\! A^R}\to {_L A_L}$ where
$R$ is a ring isomorphic to $L^{op}$ and the $R$$-$$R$-bimodule
structure of $A$ is given by fixing an isomorphism $\mu: L^{op}\to R$:
\be r\cdot a \cdot r^{\prime}\colon = a s_L\ci \theta_L^{-1}\ci \mu^{-1}(r)
t_L\ci \mu^{-1}(r^{\prime}).\lb{rbim}\ee 
The usage of the same notation ${^RA^R}$ as in (\ref{erbim}) is not
accidental. It will turn out from the next Proposition \ref{defseq}
that there exists a right bialgebroid structure on the total ring $A$
over the base $R$ for which the $R$-$R$-bimodule (\ref{erbim}) is
(\ref{rbim}). 

It makes sense to introduce the maps
\bea S_{A\ot_L A}: \atl\ot \asl &\to& \asr\ot \atr \nn
          a\ot b &\mapsto& S(b)\ot S(a) \quad {\rm and}\nn
      S_{A\ot_R A}:\asr\ot \atr &\to& \atl\ot \asl \nn
           a\ot b &\mapsto& S(b)\ot S(a).
\lb{stw}\eea
It is useful to give some alternative forms of the Definition
\ref{def}:

\bp \lb{defseq} The following are equivalent:

i) $({\cal A}_L, S)$ is a Hopf algebroid

ii) ${\cal A}_L=(A,L,s_L,t_L,\gamma_L,\pi_L)$ is a left bialgebroid
and $S$ is an anti-automorphism of the total ring $A$ satisfying
(\ref{axi}), (\ref{luap}) and
\bea 
 S_{A\ot_L A}\ci \gamma_L\ci S^{-1}&=&S_{A\ot_R A}^{-1}\ci\gamma_L \ci 
S\lb{luiii}\\ 
 (\gamma_L \ot \id_A)\ci \gamma_R = (\id_A \ot \gamma_R)\ci
\gamma_L \quad &,&\quad
 (\gamma_R \ot \id_A)\ci \gamma_L = (\id_A \ot \gamma_L)\ci
\gamma_R
\lb{luiv}
\eea
where 
we introduced the ring
$R$ and  the $R$-$R$ bimodule ${^RA^R}$ as in 
(\ref{rbim}) and the map $\gamma_R\colon =
S_{A\ot_L A}\ci {\gamma_L} \ci S^{-1} \equiv   S_{A \ot_R
A}^{-1}\ci {\gamma_L}   \ci S:A\to \asr\ot \atr$. The equations in (\ref{luiv})
are equalities of maps  $A\to \atl
\ot {\asl^R}  
\ot {\atr} $ and $A\to \asr 
\ot {\atr_L} 
\ot \asl$, 
respectively. 

iii) ${\cal A}_L=(A,L,s_L,t_L,\gamma_L,\pi_L)$ is a left bialgebroid and 
${\cal A}_R=(A,R,s_R,t_R,\gamma_R,\pi_R)$ is a right bialgebroid such
that the base 
rings are related to each other via $R\simeq L^{op}$. $S$ is a
bijection of additive groups  and 
\bea  s_L(L)=t_R(R) \quad &,&\quad t_L(L)=s_R(R) 
\quad as\ subrings\ of \ A\lb{defi}\\
 (\gamma_L \ot \id_A)\ci \gamma_R = (\id_A \ot \gamma_R)\ci
\gamma_L \quad &,&\quad
 (\gamma_R \ot \id_A)\ci \gamma_L = (\id_A \ot \gamma_L)\ci
\gamma_R \lb{defii}\\
  S(t_L(l)a t_L(l^{\prime}))=s_L(l^{\prime})S(a) s_L(l) \quad &,&\quad
            S(t_R(r^{\prime})a t_R(r))=s_R(r)S(a) s_R(r^{\prime}) 
\lb{defiii}\\
\ S(a\di)a\dii =s_R\ci \pi_R (a) \quad &,&\quad
         a\ui S(a\uii)=s_L\ci \pi_L (a) \lb{defiv}
\eea
hold true for all $l,l^{\prime}\in L$, $r,r^{\prime}\in R$ and $a\in A$.

iv) ${\cal A}_L$ is a left bialgebroid over $L$ and
${\cal A}_R$ is a right bialgebroid over $R$ such that the base
rings are related to each other via $R\simeq L^{op}$ and the
equations (\ref{defi}) and (\ref{defii}) hold true. Furthermore the
maps of additive groups
\bea \alpha: \asr \ot {_R\!A} &\to \  \atl \ot \asl \quad,\quad
a\ot b &\mapsto\  a\di \ot a\dii b \qquad and \nn
\beta: {A_R}\ot \atr &\to\ \atl \ot \asl \quad,\quad
a\ot b &\mapsto\  b\di a \ot b\dii 
\lb{schrell}
\eea
are bijective. (All modules appearing in (\ref{schrell}) are the
canonical modules introduced in Section \ref{preli}.)
\ep

Each characterization of Hopf algebroids in Proposition \ref{defseq}
will be relevant in what follows. The one in {\em ii)} is
as similar to \cite{Lu} as possible which will
be useful in Subsection \ref{ex} both in checking that concrete
examples of Lu-Hopf algebroids satisfy the axioms in Definition
\ref{def} and also in constructing Hopf algebroids in the sense of
Definition \ref{def} which do not satisfy the Lu axioms.

As it will turn out from the following proof of Proposition
\ref{defseq},
 the  Definition \ref{def} implies the existence of a
right bialgebroid structure on the total ring of the Hopf
algebroid. The characterization in {\em iii)} uses the left and right
bialgebroids underlying a Hopf algebroid in a perfectly symmetric
way. This characterization will be appropriate in developing the
theory of integrals in Section \ref{int}.

The characterization in {\em iv)} is formulated in the spirit of
\cite{Sch}, that is the bijectivity of certain Galois maps is
required. The relevance of this form of the definition is that it
shows that the Hopf algebroid in the sense of Definition \ref{def} is
a special case of Schauenburg's $\times_L$-Hopf algebra.

\pr $ ii)\Rightarrow iii):$ We construct a right bialgebroid
${\cal A}_R$ such that  $({\cal A}_L,{\cal A}_R,S)$
satisfies the  requirements in {\em iii)}: Let $R$ be a ring
isomorphic to $L^{op}$ and $\mu:L^{op}\to R$ a fixed isomorphism. Set
\be {\cal A}_R=(A,R,s_R\colon =t_L\ci \mu^{-1},t_R\colon = S^{-1}\ci
t_L\ci \mu^{-1},\gamma_R\colon = S_{A\ot_R A}^{-1} \ci {\gamma_L}
\ci S, \pi^R \colon =\mu\ci \pi_L\ci S). \lb{rbgdmu}\ee

\smallskip

$ iii)\Rightarrow iv)$
We construct the inverses of the maps (\ref{schrell}):
\bea \alpha^{-1}:  \atl \ot \asl &\to \  \asr \ot {_R A} \quad,\quad
a\ot b &\mapsto\ a\ui \ot S(a\uii)b  \qquad and \nn
\beta^{-1} :  \atl \ot \asl &\to\ {A_R}\ot \atr  \quad,\quad
a\ot b &\mapsto \    S^{-1}(b\ui)a \ot b\uii.
\lb{schinv}
\eea

\smallskip

$ iv)\Rightarrow i)$
Recall that the requirements in {\em iv)} imply that both left bialgebroids
${\cal A}_L$ and ${\cal A}_{Lcop}$ are $\times_L$-Hopf algebras in the
sense of \cite{Sch}. In particular the Proposition 3.7. of \cite{Sch}
holds true for both. That is, denoting $\alpha^{-1}(a\ot 1_A)\colon = a_+
\ot a_-$ and $\beta^{-1}(1_A \ot a) \colon = {a_{[-]}} \ot {a_{[+]}}$
we have 
\be \begin{array}{lll}
i)&{a_{+(1)}} \ot {a_{+(2)}} a_- = a\ot 1_A \quad \qquad&
{a_{[+](1)}} {a_{[-]}} \ot {a_{[+](2)}}= 1_A \ot a\nn
ii)&  {a_{(1)+}} \ot {a_{(1)-}} a\dii = a\ot 1_A \quad \qquad&
{a_{(2)[-]}} a\di \ot {a_{(2)[+]}}= 1_A \ot a \nn
iii)& (ab)_+ \ot (ab)_- = a_+ b_+ \ot b_- a_-  \quad \qquad&
{(ab)_{[-]}} \ot {(ab)_{[+]}}= {b_{[-]}}
{a_{[-]}}\ot  {a_{[+]}} {b_{[+]}}\nn 
iv)& (1_A)_+ \ot (1_A)_- = 1_A\ot 1_A \quad \qquad&
{(1_A)_{[-]}} \ot {(1_A)_{[+]}} = 1_A \ot 1_A \nn
v)& {a_{+(1)}} \ot {a_{+(2)}} \ot a_- = 
\quad \qquad&
{a_{[-]}}\ot {a_{[+](1)}} \ot {a_{[+](2)}} = 
\nn
&\qquad \qquad \qquad= a\di \ot a_{(2)+} \ot a_{(2)-}\quad 
\qquad&\qquad\qquad \qquad={a_{(1)[-]}} \ot {a_{(1)[+]}} \ot a\dii \nn
vi)& a_+ \ot {a_{-(1)}}  \ot {a_{-(2)}} = 
\quad \qquad&
 {a_{[-](2)}}\ot {a_{[-](1)}}\ot {a_{[+]}} =\nn
&\qquad\qquad\qquad = a_{++} \ot a_- \ot a_{+-}\quad \qquad&\qquad\qquad\qquad
 {a_{[-]}}\ot {a_{[+][-]}}\ot{a_{[+][+]}}\nn
vii)& a= a_+ t_L\ci \pi_L(a_-) \quad \qquad&
a= a_{[+]} s_L \ci \pi_L (a_{[-]}) \nn
viii)&  a_+ a_- = s_L \ci \pi_L (a) \quad \qquad&
a_{[+]} a_{[-]} = t_L\ci \pi_L(a).
\end{array} \lb{schpr}
\ee
 We define the antipode as 
\be S(a)\colon = s_R\ci \pi_R (a_+) a_- \lb{schantip} \ee
and what is going to be its inverse as
\be  
S^{\prime}(a)\colon = t_R\ci \pi_R (a_{[+]}) a_{[-]}. \lb{schantipinv} 
\ee
The maps (\ref{schantip}) and (\ref{schantipinv}) are  well defined due to 
the $R$-module map property of $\pi_R$.

Since $\alpha(1_A\ot s_L(l))=1_A\ot s_L(l) \equiv t_L(l) \ot 1_A$, the
requirement (\ref{axi}) holds true. By making use of {\em vi)} and
{\em i)} of (\ref{schpr}) one verifies
\bea  S(a\di)_{(1)^{\prime}}a\dii \ot  S(a\di)_{(2)^{\prime}}&=&
a_{(1)-(1)^{\prime}}a\dii\ot s_R\ci
\pi_R(a_{(1)+})a_{(1)-(2)^{\prime}}=
a_{(1)-}a\dii\ot S(a_{(1)+})=\nn
&=& 1\ot S(a)\nonumber\eea
and similarly 
\[ S^{\prime}(a\dii)_{(1)^{\prime}}\ot S^{\prime}(a\dii)_{(2)^{\prime}}
    a\di=S^{\prime}(a)\ot 1_A\]
which becomes the requirement (\ref{axsi}) once we proved
$S^{\prime}=S^{-1}$. As a matter of fact by {\em vi)} of (\ref{schpr})
\[S(a)\di \ot S(a)\dii=a_{-(1)}\ot s_R\ci \pi_R(a_+)a_{-(2)}=a_-\ot
S(a_+)\]
hence using {\em ii)} of (\ref{schpr})
\[ \beta(a\dii\ot S(a\di))=a_{(1)-}a\dii\ot S(a_{(1)+}=1_A\ot S(a),\]
so by {\em viii)} of (\ref{schpr}) and (\ref{defi})
\bea S^{\prime}\ci S(a)&=& t_R\ci \pi_R\left(s_R\ci
\pi_R(a_{(1)+})a_{(1)-}\right) a\dii=
 t_R\ci \pi_R\left(a_{(1)+}a_{(1)-}\right) a\dii=\nn
&=&s_L\ci \pi_L(a\di)a\dii=a.
\nonumber\eea
In a similar way one checks that $S\ci S^{\prime}=\id_A$. 

The anti-multiplicativity of $S$ is proven as follows:
We have $\alpha(t_R(r)\ot 1_A)= t_R(r)\ot 1_A$ and by $\beta(1_A\ot
s_R(r))= 1_A\ot s_R(r)$ and $S^{\prime}=S^{-1}$ also
$S(t_R(r)a)=S(a)s_R(r)$  hence
\bea S(ab) &=& s_R\ci \pi_R \left( (ab)_+ \right) (ab)_-=
s_R\ci \pi_R ( a_+ b_+ ) b_- a_- =
s_R\ci \pi_R \left( t_R\ci \pi_R (a_+) b_+ \right) b_- a_- =\nn
&=& s_R\ci \pi_R \left( [t_R\ci \pi_R (a_+) b]_+ \right) [t_R\ci \pi_R 
(a_+) b]_- a_- = S\left(t_R\ci \pi_R (a_+) b \right) a_- =\nn
&=& S(b) s_R\ci \pi_R ( a_+)a_- = S(b) S(a). \nonumber
\eea

\smallskip

$ i)\Rightarrow ii)$
The requirements (\ref{axi}) and (\ref{luap}) hold obviously true.
One easily checks that the maps $\alpha$ and $\beta$ in
(\ref{schrell}) are bijective with inverses
\bea \alpha^{-1}(a\ot b)&=&S^{-1}\left( S(a)\dii\right)\ot S(a)\di
b\nn
     \beta^{-1}(a\ot b) &=&S^{-1}(b)\dii a \ot S\left(
     S^{-1}(b)\di\right).
\nonumber\eea
This implies that the Proposition 3.7 in \cite{Sch} holds true both in
${\cal A}_L$ and ${\cal A}_{Lcop}$. In particular introducing the maps 
\bea \gamma_R: \ A&\to\ \asr\ot \atr\quad,\quad a&\mapsto\ (\id_A\ot
S^{-1})\ci \alpha^{-1}(a\ot 1_A)\nn
      \gamma_R^{\prime}: \ A&\to\ \asr\ot \atr\quad,\quad a&\mapsto\
      (S\ot \id_A)\ci \beta^{-1}(1_A\ot a)
\nonumber\eea
the part {\em v)} of (\ref{schpr}) reads as
\bea (\gamma_L\ot \id_A)\ci \gamma_R&=&(\id_A\ot \gamma_R)\ci
\gamma_L\lb{lem1}\\
   (\id_A\ot \gamma_L)\ci \gamma_R^{\prime}&=&(\gamma_R^{\prime}\ot
   \id_A)\ci \gamma_L. \lb{lem3}
\eea
This means that both (\ref{luiii}) and (\ref{luiv}) follow provided
$\gamma_R=\gamma_R^{\prime}$. Using the Sweedler's notation
$\gamma_R(a)=a\ui\ot a\uii$ and $\gamma_R^{\prime}(a)=a^{<1>}\ot
a^{<2>}$ by the repeated use of (\ref{lem1}) and (\ref{lem3}) we
obtain 
\bea (\id_A\ot \gamma_R)\ci \gamma_R^{\prime}(a)&=&
a^{<1>}\ot s_L\ci \pi_L({a^{<2>}}\di){{a^{<2>}}\dii}\ui \ot
{{a^{<2>}}\dii}\uii=\nn
&=&{a\di}^{<1>}\ot s_L\ci \pi_L({a\di}^{<2>}) {a\dii}\ui \ot {a\dii}\uii= \nn
&=&{{a\ui}\di}^{<1>}\ot s_L\ci \pi_L({{a\ui}\di}^{<2>}){a\ui}\dii \ot
a\uii=\nn
&=&{a\ui}^{<1>}\ot s_L\ci \pi_L({{a\ui}^{<2>}}\di){{a\ui}^{<2>}}\dii\ot
a\uii=\nn
&=&(\id_A\ot \gamma_R^{\prime})\ci \gamma_R(a).
\nonumber\eea
Since by {\em viii)} of (\ref{schpr}) both $a\ui
S(a\uii)=s_L\ci\pi_L(a)$ and $a^{<1>}S(a^{<2>})=s_L\ci\pi_L(a)$ we
have $\varepsilon(a)\colon= S\ci s_L\ci \pi_L\ci
S^{-1}(a)=S(a\di)a\dii=S^{-1}\ci s_L\ci  \pi_L\ci S(a)$ and 
\bea (m_A\ot \id_A)&\ci&(\id_A\ot \varepsilon\ot \id_A)\ci(\id_A\ot
\gamma_R)\ci  \gamma_R^{\prime}(a)=
a^{<1>}\varepsilon({a^{<2>}}\ui)\ot {a^{<2>}}\uii=\nn
&=&a^{<1>}\ot {a^{<2>}}\uii  S^{-2}\ci s_L\ci  \pi_L\ci S({a^{<2>}}\ui)=\nn
&=&a^{<1>}\ot S^{-1}\left( S(a^{<2>})\di\right) S^{-1}\ci t_L\ci
\pi_L\left( S(a^{<2>})\dii\right)=  \gamma_R^{\prime}(a)\nn
(m_A\ot \id_A)&\ci&(\id_A\ot \varepsilon\ot
\id_A)\ci(\gamma_R^{\prime}\ot \id_A)\ci \gamma_R(a)=
{a\ui}^{<1>}\varepsilon({a\ui}^{<2>})\ot a\uii=\nn
&=&{a\ui}^{<1>} S\ci s_L\ci \pi_L\ci S^{-1}({a\ui}^{<2>})\ot a\uii=\nn
&=&S\left(S^{-1}(a\ui)\dii\right)S\ci s_L\ci
\pi_L\left(S^{-1}(a\ui)\di\right)\ot a\uii= 
\gamma_R(a)\nonumber\eea
hence by the equality of the left hand sides
$\gamma_R=\gamma_R^{\prime}$.
\hspace{1cm}\qed

The following is a consequence of the proof of Proposition
\ref{defseq}:

\bp \lb{sbgdi} Let $({\cal A}_L,S)$ be a Hopf algebroid and ${\cal
A}_R$ a right bialgebroid such that $({\cal A}_L,{\cal
A}_R,S)$ satisfies the requirements in {\em iii)} of Proposition
\ref{defseq}. Then
both $\left(S:A\to A^{op},\nu\colon =\pi_R \ci s_L:\right.$\break 
$\left. L\to R\op\right)$ and 
$\left(S^{-1}:A\to A\op,\mu\colon = \pi_R\ci t_L:L\to 
R^{op}\right)$ are left bial\-ge\-broid
iso\-mor\-phisms ${\cal A}_L \to ({\cal A}_R)_{cop}^{op}$.  In particular
${\cal A}_R$ is unique up to an isomorphism of the form $(\id_A,\phi)$.
\ep

One easily checks that $\mu^{-1} \ci \nu = \theta_L$. For the sake of
symmetry we introduce also $\theta_R\colon = \nu\ci \mu^{-1}$ with the
help of which the right analogue of (\ref{theta}) holds true:
\[ S\ci s_R= t_R \ci \theta_R. \]

Proposition \ref{sbgdi} has an interpretation in terms of the
forgetful functors $\Phi_R :{\cal M}_A \to {_R {\cal M}_R}$ and 
$\Phi_L :{_A{\cal M}} \to {_R {\cal M}_R}$ as follows.
The antipode map defines two  functors ${\cal S}$ and ${\cal
S}^{\prime}: {\cal M}_A \to {_A {\cal M}}$. They have object maps
$(M, \triangleleft)\mapsto (M,\triangleleft \ci S)$ and 
$(M,\triangleleft)\mapsto (M,\triangleleft \ci S^{-1})$, respectively,
and the identity maps on the morphisms. It is 
clear that ${\cal S}$  and  ${\cal S}^{\prime}$ are  strict antimonoidal
equivalence functors. 
The ring automorphisms $\mu$ and $\nu$ define endo-functors ${\underline \mu}$ 
and ${\underline \nu}$ of ${_R {\cal M}_R}$. 
The object maps are $(M, \triangleright, \triangleleft)
\mapsto (M, \triangleleft\ci \mu, \triangleright\ci\mu )$ 
and  $(M, \triangleright, \triangleleft)\mapsto 
(M, \triangleleft\ci \nu, \triangleright\ci\nu )$, respectively, and
the identity map on the morphisms. The ${\underline \mu}$ and
${\underline \nu}$ are also strict antimonoidal equivalence functors. 
We have then equalities of strong monoidal functors: $\Phi_L \ci {\cal
S} = {\underline \nu} \ci \Phi_R$ 
and $\Phi_L \ci {\cal S^{\prime}} = {\underline \mu} \ci \Phi_R$.

Finally we define the morphisms of Hopf algebroids:
\bd \lb{hgdeq}
A {\em Hopf algebroid homomorphism (isomorphism)}  $({\cal
A}_L,S)\to ({\cal A}_L^{\prime},S^{\prime})$ is a left bialgebroid
homomorphism (isomorphism) ${\cal A}_L\to {\cal A}_L^{\prime}$. A Hopf
algebroid homomorphism $(\Phi,\phi)$ is {\em strict} if $S^{\prime}\ci
\Phi=\Phi\ci S$. 
\ed
The existence of non-strict isomorphisms of Hopf algebroids -- that is
the non-uniqueness of the antipode in a Hopf algebroid -- is a new
feature compared to (weak) Hopf algebras. The antipodes making a given
left bialgebroid into a Hopf algebroid are characterized in \cite{B}.

In the following (in particular in Section \ref{int}) we are going to
call a triple $({\cal A}_L,{\cal A}_R,S)$ satisfying the {\em iii)} of
Proposition \ref{defseq} a {\em symmetrized form} of the Hopf
algebroid $({\cal A}_L,S)$. The ${\cal A}_R$ is called the right
bialgebroid underlying $({\cal A}_L,S)$. In the view of Proposition
\ref{sbgdi} the symmetrized form is unique up to the choice of the base
ring $R$ of ${\cal A}_R$ within the isomorphism class of $L^{op}$ -- the
opposite of the base ring of ${\cal A}_L$ --  and the isomorphism
$\mu:L^{op}\to R$ in (\ref{rbgdmu}). 

Let us define the {\em homomorphisms} of symmetrized Hopf algebroids
${\cal A}=({\cal A}_L,{\cal A}_R,S)\to {\cal A}^{\prime}=({\cal
A}_L^{\prime}, {\cal A}_R^{\prime},S^{\prime})$ as pairs of
bialgebroid homomorphisms $(\Phi_L,\phi_L): {\cal A}_L\to {\cal
A}_L^{\prime}$,  $(\Phi_R,\phi_R): {\cal A}_R\to {\cal
A}_R^{\prime}$. Then by Proposition \ref{sbgdi} the Hopf algebroid
homomorphisms $(\Phi,\phi):({\cal A}_L,S)\to ({\cal A}_L^{\prime},
S^{\prime})$ are injected into the homomorphisms of symmetrized Hopf
algebroids $({\cal A}_L,{\cal A}_R,S)\to ({\cal
A}_L^{\prime}, {\cal A}_R^{\prime},S^{\prime})$ via
$(\Phi,\phi)\mapsto \left( (\Phi,\phi)\ ,\ ((S^{\prime-1}\ci \Phi\ci
S,\mu^{\prime}\ci \phi\ci \mu^{-1})\right)$ -- where $\mu=\pi_R\ci
t_L$ and $\mu^{\prime}=\pi_R^{\prime}\ci t^{\prime}_L$ are the ring
isomorphisms introduced in Proposition \ref{sbgdi}. This implies that
two symmetrized Hopf algebroids $({\cal A}_L,{\cal A}_R,S)$ and $({\cal
A}_L^{\prime}, {\cal A}_R^{\prime},S^{\prime})$ are isomorphic if and
only if the Hopf algebroids $({\cal A}_L,S)$ and $({\cal
A}_L^{\prime}, S^{\prime})$ are isomorphic.

A homomorphism $\left( (\Phi_L,\phi_L)\ ,\
(\Phi_R,\phi_R)\right)$ of symmetrized Hopf algebroids ${\cal A}\to
{\cal A}^{\prime}$  is {\em strict} if $\Phi_L=\Phi_R$ as
homomorphisms of rings $A\to A^{\prime}$. We leave it to the reader to
check 
that this is equivalent to the requirement that $ (\Phi_L,\phi_L)$ is
a strict homomorphism of Hopf algebroids $({\cal A}_L,S)\to ({\cal
A}_L^{\prime}, S^{\prime})$, that is two symmetrized Hopf algebroids
$({\cal A}_L,{\cal A}_R,S)$ and $({\cal 
A}_L^{\prime}, {\cal A}_R^{\prime},S^{\prime})$ are strictly isomorphic
if and 
only if the Hopf algebroids $({\cal A}_L,S)$ and $({\cal
A}_L^{\prime}, S^{\prime})$ are strictly isomorphic.

The usage of symmetrized Hopf algebroids allows for the definition of
the opposite and co-opposite structures: ${\cal A}^{op}=({\cal A}_R^{op},
{\cal A}_L^{op},S^{-1})$  and  ${\cal A}_{cop}=({\cal A}_{L\ cop},
{\cal A}_{R\ cop},S^{-1})$, respectively.

\subsection{Relation to the Hopf bialgebroid of Day and Street}
\lb{DS}

In \cite{Sz: EM} left bialgebroids over the base $L$ have been characterized as
opmonoidal monads $\mathtt{T}$ on the category $_L\M_L$ such that their
underlying functors $T$ have right adjoints. The endofunctor $T$ is given by
tensoring over $L^e=L\ot L\op$ with the $L^e$-$L^e$-bimodule $A$. The monad
structure makes $A$ into an $L^e$-ring via an algebra map $\eta\colon L^e\to
A$, while the opmonoidal structure comprises the coproduct and counit of the
bialgebroid $A$. In a recent preprint \cite{DS} Day and Street put this into
the more general context of pseudomonoids in monoidal bicategories. Using the
notion of $*$-autonomy \cite{Barr} they propose a definition of
{\em Hopf bialgebroid} as a $*$-autonomous structure on a bialgebroid. We are
going to show in this section that their definition coincides with our
Definition \ref{def}.

Working with $k$-algebras the base category $\V$ is the category $\M_k$ of
$k$-modules. Then the monoidal bicategory in question is $\Mod(\V)$ having as
objects the $k$-algebras and as hom-categories $\Mod(\V)(A,B)$ the category
$_B\M_A$ of $B$-$A$-bimodules. The tensor product $\ot=\ot_k$ of $k$-modules
makes $\Mod(\V)$ monoidal. A pseudomonoid in $\Mod(\V)$ is a triple $\bra
A,M,J\ket$ where $A$ is an object, $M\colon A\ot A\to A$ and $J\colon k\to A$
are 1-cells satisfying associativity and unitality up to invertible 2-cells
that in turn satisfy the pentagon and the triangle constraints. If $\bra
A,L,s,t,\gamma,\pi\ket$ is a left bialgebroid in $\V$ then we have a
strong monoidal morphism 
\begin{equation}
\eta^*\colon\bra A,M,J\ket\ \to\ \bra L^e,m,j\ket
\end{equation}
of pseudomonoids where
\begin{eqnarray}
M&=&A\ot_L  A\ \text{with actions}\ a\cdot(x\ot_L  y)\cdot(a_1\ot a_2)=
a\di xa_1\ot_L a\dii ya_2\\
J&=&L\ \text{with action}\ a\cdot x=\pi(as(x))\\
m&=&L^e\ot_L  L^e\ \text{with actions}\ (l\ot l')\cdot((x\ot x')\ot_L (y\ot y'))
\cdot((l_1\ot l'_1)\ot(l_2\ot l'_2))\nonumber\\
&&\ \ \ \ \ \ \ \ \ \ \ \ \ \ \ \ \ \ \ \ \ \ \ \ \ \ \ = (lxl_1\ot
l'_1x')\ot_L  (yl_2\ot l'_2y'l')\\ 
j&=&L\ \text{with action}\ (l\ot l')\cdot x=lxl'
\end{eqnarray}
and the bimodule $\eta^*=\,_{L^e}A_A$, induced by the algebra map $\eta$, has
a left adjoint $\eta_*=\,_AA_{L^e}$ with the counit $\epsilon\colon
\eta_*\circ\eta^*\to A$ induced by multiplication of $A$. The connection of
this bimodule picture with module catgeories can be seen by applying the
monoidal pseudofunctor  $\Mod(\V)(k,-)\colon\Mod(\V)\to\Cat$. Then
pseudomonoids become monoidal categories, $\eta^*$ becomes the strong monoidal
forgetful functor $\eta^*\ot_A-\colon_A\M\to\,_L\M_L$ and $\eta_*$ its
opmonoidal left adjoint. An important piece of the structure is the bimodule
morphism
\begin{eqnarray} \lb{psi map}
\psi&\colon&\eta_*\circ m\to M\circ(\eta_*\ot\eta_*)\\
a\ot(l\ot l')&\mapsto& a\di s(l)\ot a\dii t(l')\nonumber
\end{eqnarray}
which makes $\eta_*$ opmonoidal and encodes the comultiplication $\gamma\colon
A\to A\ot_L A$ of the bialgebroid.

In \cite{DS} a Hopf bialgebroid is defined to be a bialgebroid $A$ over $L$
together with a strong $*$-autonomous structure on the
opmonoidal morphism $\eta_*\colon L^e\to A$. The latter means
\begin{enumerate}
\item a $*$-autonomous structure on the pseudomonoid $A$, i.e., 
\begin{enumerate}
\item a right $A\ot A$-module $\sigma$ defined on the
$k$-module $A$ in terms of an algebra isomorphism $\xi\colon A\to A\op$ by
\begin{equation} \lb{sigma module}
x\cdot(a\ot b)\ =\ \xi^{-1}(b)x a\,,\qquad x,a,b\in A
\end{equation}
\item and an isomorphism 
\[
\Gamma\colon \sigma\circ(M\ot A)\to \sigma\circ(A\ot M)
\]
\end{enumerate}
\item a "canonical" $*$-autonomous structure on the pseudomonoid $L^e$ which
consists of
\begin{enumerate}
\item a right $L^e\ot L^e$-module $\sigma_0$ defined on $L^e$ using the
isomorphism $\xi_0\colon L^e\to (L^e)\op$, $l\ot l'\mapsto\theta^{-1}(l')\ot l$
for some algebra automorphism $\theta\colon L\to L$. Thus
\[
(x\ot x')\cdot((l_1\ot l_1')\ot(l_2\ot l_2'))\ =\ \theta^{-1}(l_2')xl_1\ot
l_1'x'l_2
\]
for $x,x',l_1,l_1',l_2,l_2'\in L$ where juxtaposition is always multiplication
in $L$ and never in $L\op$. (Note that allowing a non-trivial $\theta$ in
the definition of $\sigma_0$ is a slight deviation from \cite{DS} which is,
however, well motivated by Section \ref{ringex}.)
\item and the isomorphism
\[
\Gamma_0\colon\sigma_0\circ(m\ot L^e)\to\sigma_0\circ(L^e\ot m)
\]
\end{enumerate}
\item and the arrow $\eta_*$ is strongly $*$-autonomous in the
following sense. There exists a 2-cell $\tau\colon
\sigma_0\circ(\eta_*\ot\eta_*)\to\sigma$ such that
\begin{equation} \lb{*-auto hexagon}
[\Gamma_0\circ(\eta_*\ot\eta_*\ot\eta_*)]\bullet[\sigma\circ(\psi\ot\eta_*)]
\bullet[\tau\circ(m\ot L^e)]
=[\sigma\circ(\eta_*\ot\psi)]\bullet[\tau\circ(L^e\ot m)]\bullet\Gamma_0
\end{equation}
and such that 
\begin{eqnarray} \lb{tau mate}
\tau^l&\colon& \sigma_0\circ(\eta^*\ot L^e)\to \sigma\circ(A\ot\eta_*)\\
\tau^l&=&[\sigma\circ(\epsilon\ot\eta_*)]\bullet[\tau\circ(\eta^*\ot L^e)]
\nonumber
\end{eqnarray}
is an isomorphism. (We used $\circ$ and $\bullet$ to denote horizontal and
vertical compositions in the bicategory $\Mod(\V)$.) 
\end{enumerate}

Now we are going to translate this categorical definition into simple algebraic
expressions. 
\bl \lb{*autA}
$^*$-autonomous structures on the pseudomonoid $\bra A,M,J\ket$ are in
one-to-one correspondence with the data $\bra\xi, \sum_k e_k\ot f_k\ket$ 
where $\xi$ is an anti-automorphism of the ring $A$ and $\sum_k e_k\ot f_k
\in A \stackrel{\ot}{_{_L}} A$ is such that for all $a\in A$
\begin{equation} 
\sum_k\xi(a\di)_{(1^{\prime})} e_k a\dii \ot
\xi(a\di)_{(2^{\prime})} f_k \ =\ \sum_k e_k\ot f_k  \xi(a)\lb{preaxsi}
\end{equation}
and such that there exists $\sum_j g_j\ot h_j\in A\stackrel{\ot}{_{_L}}A$
satisfying for all $a\in A$
\bea
\sum_j\xi^{-1}(a\dii)_{(1^{\prime})} g_j \ot \xi^{-1}(a\dii)_{(2^{\prime})}
h_j a\di &=& \sum_j g_j 
\xi^{-1}(a) \ot h_j \lb{preaxs}\\
\sum_{j,k}\xi(g_j)\di e_k h_j \ot \xi(g_j)\dii f_k &=& 1_A\ot 1_A \lb{ggi}\\
\sum_{j,k}\xi^{-1}(f_k)\di g_j \ot \xi^{-1}(f_k)\dii h_j e_k &=& 1_A\ot 1_A 
\lb{gig}
\eea
as elements of $A\stackrel{\ot}{_{_L}}A$.  
\el
\pr In order to find explicit formulas for $\Gamma$ and its inverse we
introduce the isomorphisms 
\begin{eqnarray}
\varphi_+\colon\ \sigma\circ(M\ot A)&\to& A\stackrel{\ot}{_{_L}}A\\
x\stackrel{\ot}{_{_{A\ot A}}}(a\stackrel{\ot}{_{_L}}b\ot c)&\mapsto&
(\xi^{-1}(c)x)\di a\stackrel{\ot}{_{_L}}(\xi^{-1}(c)x)\dii b\nn\\
\varphi_-\colon\ \sigma\circ(A\ot M)&\to&A\stackrel{\ot}{_{_L}}A\\
x\stackrel{\ot}{_{_{A\ot A}}}(a\ot b\stackrel{\ot}{_{_L}}c)&\mapsto&
\xi(xa)\di b\stackrel{\ot}{_{_L}}\xi(xa)\dii c\nonumber
\end{eqnarray}
with inverses
\begin{eqnarray*}
\varphi_+^{-1}(a\stackrel{\ot}{_{_L}}b)&=&1\stackrel{\ot}{_{_{A\ot
A}}}(a\stackrel{\ot}{_{_L}}b\ot 1)\\
\varphi_-^{-1}(b\stackrel{\ot}{_{_L}}c)&=&1\stackrel{\ot}{_{_{A\ot
A}}}(1\ot b\stackrel{\ot}{_{_L}}c)
\end{eqnarray*}
Then $\Gamma':=\varphi_-\circ\Gamma\circ\varphi_+^{-1}$ is a twisted $A^{\ot
3}$-automorphism of $A\ot_L A$ in the sense of
\begin{equation}
\Gamma'(\xi^{-1}(c)\cdot(x\stackrel{\ot}{_{_L}}y)\cdot(a\ot b))
=\xi(a)\cdot \Gamma'(x\stackrel{\ot}{_{_L}}y)\cdot (b\ot c)
\end{equation}
Hence, $\Gamma$ is uniquely determined by $\sum_i e_i\ot
f_i=\Gamma'(1\ot 1)$ as
\begin{equation} \lb{Gamma map}
\Gamma(x\stackrel{\ot}{_{_{A\ot A}}}(a\stackrel{\ot}{_{_L}}b\ot c))
=1\stackrel{\ot}{_{_{A\ot A}}}(x\di a\ot\sum_k e_k x\dii
b\stackrel{\ot}{_{_L}}f_k c)
\end{equation}
and $\sum_k e_k\ot f_k$ satisfies (\ref{preaxsi}). Invertibility then
implies that ${\Gamma'}^{-1}(1\ot 1)=\sum_j g_j\ot h_j$ satisfies the
remaining equations. 
\qed

We need also the expression for $\Gamma_0$. We leave it to the reader to check
that
\begin{eqnarray} \lb{Gamma0 map}
&&\Gamma_0((x\ot x')\stackrel{\ot}{_{_{L^e\ot L^e}}}
(l_1\ot l'_1)\stackrel{\ot}{_{_L}}(l_2\ot l'_2)\ot(l_3\ot l'_3))=\\
&=&(1\ot 1)\stackrel{\ot}{_{_{L^e\ot L^e}}}
(xl_1\ot l'_1)\ot(l_2\ot l'_2x')\stackrel{\ot}{_{_L}}(l_3\ot l'_3)\nonumber
\end{eqnarray}
is well-defined and is invertible.
\bl \lb{*aut}
The $*$-autonomous property of $\eta_*$ is equivalent to the
existence of an element $i\in A$ satisfying
\bea i\eta(\theta^{-1}(l')\ot l)&=& \xi^{-1}(\eta(l\ot l'))i\,,\quad l,l'\in L
 \lb{iint}\\
\sum_k e_k\ot f_k \xi(i)&=& \xi(i)\di \ot \xi(i)\dii \lb{prel}
\eea
while strong $*$-autonomy adds the requirement that $i$ be invertible. 
\el
\pr
Since $1\ot 1$ is a cyclic vector in the $L^e\ot L^e$-module $\sigma_0$, the
module map $\tau$ is uniquely determined by $i:=\tau(1\ot 1)$ as
\begin{equation} \lb{tau map}
\tau(x\ot x')\ =\ i\eta(x\ot x')\,.
\end{equation}
Tensoring with $\eta_*=\,_AA_{L^e}$ from the right being the restriction via
$\eta\colon L^e\to A$ the element $i$ is subject to condition (\ref{iint}) due
to (\ref{sigma module}). Using the expressions (\ref{tau map}), (\ref{Gamma
map}), (\ref{Gamma0 map}) and (\ref{psi map}) the $*$-autonomy condition
(\ref{*-auto hexagon}) becomes equation (\ref{prel}). The mate of $\tau$ in
(\ref{tau mate}) can now be written as 
\[
\tau^l((x\ot x')\stackrel{\ot}{_{_{L^e\ot L^e}}}(a\ot(l\ot l')))
=i\eta(x\ot x')a\stackrel{\ot}{_{_{A\ot L^e}}}(1\ot(l\ot l'))
\]
Hence $\tau^l$ is invertible iff the map $a\mapsto ia$ is, i.e., iff $i$ is
invertible. \qed

\bt
A strong $*$-autonomous structure on the bialgebroid $A$ over $L$ in the sense
of \cite{DS} is equivalent to a Hopf algebroid structure $({\cal A}_L,S)$ in
the sense of Definition \ref{def}.
\et

\pr Using invertibility of $i\in A$ condition (\ref{iint}) has the equivalent 
form
\be \xi\left( i\eta(l'\ot l)i^{-1}\right)=\eta(l\ot\theta(l'))
\lb{dsaxi}\ee
and (\ref{prel}) can be used to express the element $\sum_k e_k\ot f_k$
in terms of $i$,
\[ e_i\ot f_i = \xi(i)\di \ot \xi(i)\dii \xi(i)^{-1}.\]
The conditions (\ref{preaxsi}-\ref{gig}) are then equivalent to
\bea g_j\ot h_j&=& i\di  i^{-1}\ot i\dii \\
\left[ i^{-1}\xi^{-1}(a\dii) i\right]_{(1)^{\prime}} \ot 
\left[ i^{-1}\xi^{-1}(a\dii) i\right]_{(2)^{\prime}}  a\di &=& 
i^{-1}\xi^{-1}(a) i \ot 1_A \lb{dsaxsi}\\
\xi\left(ia\di i^{-1}\right)_{(1)^{\prime}} 
a\dii \ot \xi\left(ia\di i^{-1}\right)_{(2)^{\prime}}  &=& 
1_A \ot \xi\left(ia i^{-1}\right).
\lb{dsaxs}
\eea
This means that introducing $S:A\to A^{op}$, $a\mapsto \xi\left(ia
i^{-1}\right)$, the conditions (\ref{dsaxi}), (\ref{dsaxsi}) and
(\ref{dsaxs}) are identical to the axioms (\ref{axi}), (\ref{axsi})
and (\ref{axs}),  respectively. 
\hspace{1cm}\qed


\subsection{Examples}
\lb{ex}

In addition to our motivating example in Section \ref{ringex} let us
collect some more examples of Hopf algebroids.
\bex \lb{whaex}
{\em Weak Hopf Algebras with bijective antipode }

Let $(H,\Delta,\varepsilon,S)$ be a weak Hopf algebra \cite{BSz,N,BNSz}
(WHA) over the commutative ring  $k$ with bijective antipode. This
means that $H$ is 
an associative unital 
$k$ algebra, $\Delta: H\to H\stackrel{\ot}{_{_k}} H$ is a coassociative
coproduct. It is an algebra map (i.e. multiplicative) but not unit
preserving in general. In its stead we have {\em weak}
comultiplicativity of the unit:
\[ 1_{[1]} \ot 1_{[2]} 1_{[1^{\prime}]} \ot 1_{[2^{\prime}]} =
1_{[1]} \ot 1_{[2]} \ot 1_{[3]}=
1_{[1]} \ot 1_{[1^{\prime}]} 1_{[2]} \ot  1_{[2^{\prime}]} \]
where $1_{[1]}\ot 1_{[2]}=\Delta(1)$.
The map $\varepsilon:H\to k$ is the  counit  of the coproduct
$\Delta$. Instead of being multiplicative it is {\em weakly}
multiplicative:  
\[ \varepsilon (ab_{[1]}) \varepsilon(b_{[2]} c)=
\varepsilon(abc)=
\varepsilon (ab_{[2]}) \varepsilon(b_{[1]} c)\qquad for\ a,b,c\in H\]
The bijective map $S:H\to H$ is the antipode, subject to the axioms
\bea h_{[1]} S(h_{[2]}) &=& \varepsilon( 1_{[1]} h)1_{[2]}\nn
S(h_{[1]}) h_{[2]} &=& 1_{[1]} \varepsilon(h1_{[2]})\nn
S(h_{[1]})h_{[2]} S(h_{[3]}) &=& S(h)
\nonumber \eea
for all $h\in H$.  
(If $H$ is finite over $k$ then the assumption made about
the bijectivity of $S$ is redundant.) The WHA
$(H,\Delta,\varepsilon,S)$ is a Hopf algebra if and only if $\Delta$ is
unit preserving.

The algebra $H$ contains two commuting subalgebras:
$R$ is the image of $H$ under the projection $\sqcap^R: h\mapsto 1_{[1]}
\varepsilon(h1_{[2]})$ and $L$ under $\sqcap^L: h\mapsto  \varepsilon(
1_{[1]} h)1_{[2]}$ -- generalizing the subalgebra of the scalars in 
a Hopf algebra. Both maps $S$ and $S^{-1}$ restrict to algebra
anti-isomorphisms $R\to L$. We have four commuting actions of $L$ and
$R$  on H: 
\bea H^R:&\quad& h\cdot r \colon = h r\nn
     {^RH}: &\quad& r\cdot h \colon = h S^{-1}(r) \nn
     {H_L}: &\quad& h\cdot l \colon = S^{-1}(l) h \nn
     {_LH}: &\quad& l\cdot h \colon = l h.
\nonumber \eea
Introduce the canonical projections $\p_R: H\stackrel{\ot}{_{_k}} H\to
H^R\ot {^RH}$ and $\p_L: H\stackrel{\ot}{_{_k}} H \to {H_L}\ot {_L H}.$ 
There exists  a left and a right bialgebroid structure corresponding
to the weak Hopf algebra:
\bea {\cal H}_R&\colon =& (H,R,\id_R,S^{-1}\vert_R, \p_R\ci\Delta,
\sqcap^R)\nn
     {\cal H}_L&\colon =&
(H,L,\id_L\,,S^{-1}\vert_L\,,\p_L\ci\Delta,\sqcap^L).
\nonumber \eea
We leave it as an exercise to the reader to check that $({\cal
H}_L,{\cal H}_R,S)$ satisfies the requirements of
Proposition \ref{defseq}.{\em iii)}.
\eex
Notice that the examples of the above class are not necessarily finite 
dimensional and not even finitely generated over $R$. To have a
trivial counterexample think of the group Hopf algebra $kG$ of an
infinite group. 

\bex \lb{twisthopfex}
{\em An example that does not satisfy the Lu-axioms \cite{Lu}}:

Let $k$ be a field the characteristic of which is different from
2. Consider the group bialgebra $kZ_2$ with presentation
\[ kZ_2\ =\ 
\text{bialg-}\langle\,t\,\vert\,t^2=1,\,
      \Delta(t)=t\ot t,\,\varepsilon(t)=1\rangle
\] 
as a left bialgebroid over the base $k$. That is to say, we set ${\cal
A}_L=(kZ_2,k,\eta,\eta,\Delta,\varepsilon)$ where $\eta$ is the unit
map $k\to kZ_2$, $\lambda\mapsto \lambda 1$. Introduce the
would-be-antipode $S:kZ_2\to kZ_2$, $t\mapsto -t$.
\bp The pair $({\cal A}_L,S)$ in the Example \ref{twisthopfex}
satisfies the axioms in Definition \ref{def} but not the Lu-axioms.
\ep

\pr One easily checks the conditions in {\em ii)} of Proposition
\ref{defseq} on the single algebraic generator $t$, proving that
$({\cal A}_L,S)$ is a Hopf algebroid in the sense of Definition
\ref{def}. 

Now, the base ring $L$ being $k$ itself, the canonical projection
$kZ_2\stackrel{\ot}{_{_{k}}} kZ_2 \to kZ_2\stackrel{\ot}{_{_{L}}}
kZ_2$ is the identity map leaving us with 
the only section $kZ_2\stackrel{\ot}{_{_{L}}} kZ_2 \to
kZ_2\stackrel{\ot}{_{_{k}}}  kZ_2$, the
identity map. Since $t\di S(t\dii)=-1$ and $\eta\ci \varepsilon (t)=1$,
this contradicts to that $({\cal A}_L,S)$ is a Lu-Hopf algebroid.
\qed
\eex
In the Example \ref{whaex} the left and right coproducts $\gamma_L$
and $\gamma_R$ are
compositions of a coproduct $\Delta:A\to A\stackrel{\ot}{_{_k}} A$
with the
canonical projections $\p_L$ and $\p_R$, respectively. Actually many
other examples can be found this way -- by allowing for $\Delta$
not to be counital:
\bp \lb{luex}
Let ${\cal A}_L$ be a left bialgebroid such that $A$ and $L$ are
$k$-algebras over some commutative ring $k$, $S$ an antiautomorphism of
$A$ such that the axioms (\ref{axi}) and (\ref{luap}) hold true.
Suppose that $\gamma_L
=\p_L\ci \Delta$, where $\p_L$ is the canonical projection
$A\stackrel{\ot}{_{_k}} A \to \atl \ot \asl$ and
$\Delta : A\to A\stackrel{\ot}{_{_k}} A$ is a coassociative (possibly
non-counital) coproduct satisfying
\bea \p_L \ci (S\ot S)\ci \Delta^{op} &=& \p_L\ci \Delta \ci
S\nn 
     \p_L \ci (S^{-1}\ot S^{-1})\ci \Delta^{op} &=& \p_L\ci \Delta \ci S^{-1}.
\lb{deltacond}\eea
Then $({\cal A}_L,S)$ is a Hopf algebroid in the sense of Definition
\ref{def}.  
\ep 
\pr We leave it to the reader to check that all the requirements of 
Proposition \ref{defseq}.{\em ii)} are satisfied.

\bex\lb{grex}{\em The groupoid Hopf algebroid }

Let ${\cal G}$ be a groupoid that is a small category with all morphisms
invertible. Denote the object set by ${\cal G}^0$ and the set of
morphisms by ${\cal G}^1$. For a commutative ring $k$ the groupoid
algebra is the 
$k$-module spanned by the elements of ${\cal G}^1$ with the
multiplication given by the composition of the morphisms if the latter
makes sense and $0$ otherwise. It is an associative algebra and if ${\cal
G}^0$ is finite it has a unit $1= \sum_{a\in {\cal G}^0} a$. The groupoid
algebra admits a left bialgebroid structure over the base subalgebra $k{\cal
G}^0$. The map $s_L=t_L$ is the canonical embedding, $\gamma_L$ is the diagonal
map $g\mapsto g\stackrel{\ot}{_{_L}} g$ and $\pi_L(g)\colon =
target(g)$. This left bialgebroid 
together with the antipode $S(g)\colon = g^{-1}$ is a Hopf algebroid
in the sense of the
Definition \ref{def}. Actually this example is of the kind described in
Proposition \ref{luex} with $\Delta(g)\colon =
g\stackrel{\ot}{_{_k}} g$. 
\eex
\bex \lb{torusex}{\em The algebraic quantum torus}

Let $k$ be a field and $T_q$ the  unital associative $k$ algebra
generated by two invertible elements $U$ and $V$ subject to the
relation $UV=qVU$ where $q$ is an invertible element in $k$.
As it is explained in \cite{KhalRang}, the algebra $T_q$ admits a
Lu-Hopf algebroid structure over the base  
subalgebra $L$ generated by $U$: the map
$s_L=t_L$ is the 
canonical embedding, $\gamma_L(U^n V^m )\colon = U^n V^m
\stackrel{\ot}{_{_L}} V^m \equiv V^m  \stackrel{\ot}{_{_L}} U^n V^m$, 
$\pi_L(U^n V^m) \colon = U^n$ and the antipode $S(U^n V^m)\colon =
V^{-m} U^n$. The section $\xi$ of the canonical projection $\p_L:T_q
\stackrel{\ot}{_{_k}} T_q \to T_q \stackrel{\ot}{_{_L}} T_q$ appearing
in the Lu axioms is of the form $\xi \left(U^nV^m \stackrel{\ot}{_{_L}} U^k
V^l \right)\colon = U^{(n+k)} V^m \stackrel{\ot}{_{_k}}  V^l$.

The reader may check that these maps satisfy the Definition \ref{def}
as well. This example is also of the type considered in 
Proposition \ref{luex} with $\Delta(U^n V^m)\colon = U^n V^m
\stackrel{\ot}{_{_k}} V^m$.
\eex
\bex \lb{BMex}{\em Examples by Brzezinski and Militaru \cite{BreMi}}

In the paper \cite{BreMi} a wide class of examples of Lu-Hopf
algebroids is described. Some other examples \cite{Lu,Pa} turn out to
belong also to this class. 

The examples of \cite{BreMi} are Lu-Hopf algebroids of the type
considered in Proposition 
\ref{luex}: let $(H,\Delta_H,\varepsilon_H,\tau)$ be a Hopf algebra over
the field $k$ with 
bijective antipode $\tau$, and the triple $(L,\cdot,\rho)$ a braided
commutative algebra in the category ${_H {\cal D}^H}$
of Yetter-Drinfel'd modules over $H$. 
Then the crossed product algebra $L\# H$ carries a left bialgebroid
structure over the base algebra  $L$: 
\bea s_L(l) &=& l\# 1_H \nn
     t_L(l) &=& \rho(l) \equiv l_{\langle 0 \rangle} \#  l_{\langle 1
            \rangle} \nn
     \gamma_L(l\#h) &=& (l\# h\di) \stackrel{\ot}{_{_L}} (1_L\# h\dii)\nn
     \pi_L(l\#h) &=& \varepsilon_H (h) l
\lb{BMbgd}\eea
where $h\di \stackrel{\ot}{_{_k}} h\dii \equiv \Delta_H(h).$ It is
proven in \cite{BreMi} that the left bialgebroid (\ref{BMbgd}) 
and the bijective antipode 
\be S(l\# h)\colon = \left( \tau (h\dii) \tau^2(l_{\langle
1\rangle})\right)\cdot l_{\langle 0\rangle} \# \tau(h\di)
\tau^2(l_{\langle 2\rangle}) \lb{BMap}\ee
form a Lu-Hopf algebroid. It is obvious that 
$\gamma_L$ is of the form $\p_L\ci \Delta$ with $\Delta(l\# h) \colon =
(l\# h\di) \stackrel{\ot}{_{_k}} (1_L\# h\dii)$. The map $\Delta$ is
well defined since $L\# H$ is $L \stackrel{\ot}{_{_k}} H$ as a k-space
and $\Delta_H$ maps $H$ into $H\stackrel{\ot}{_{_k}} H$. We leave it
to the reader to check that $\Delta $ satisfies (\ref{deltacond})
hence the left bialgebroid (\ref{BMbgd}) and the antipode (\ref{BMap})
form a Hopf algebroid in the sense of Definition \ref{def}.
\eex
The Example \ref{twisthopfex} is not of the type considered in
Proposition \ref{luex}. Although $\gamma_L$ is of the form $\p_L\ci
\Delta$, the $\Delta$ does not satisfy (\ref{deltacond}). In
\cite{KhalRang} data $({\cal A}_L,S,{\tilde S})$ satisfying 
compatibility conditions somewhat analogous to (\ref{deltacond}) were
introduced under the name {\em extended Hopf algebra}. The next
Proposition states that
extended Hopf algebras with $S$ bijective (such as Example
\ref{twisthopfex}) provide examples of Hopf algebroids:
\bp \lb{exthopf}
Let $({\cal A}_L,S,{\tilde S})$ be an extended Hopf algebra. This
means that
${\cal A}_L$ is a left bialgebroid such that $A$ and $L$ are
$k$-algebras over some commutative ring $k$. The maps  $S$ and ${\tilde
S}$ are anti-automorphisms of 
the algebra $A$, ${\tilde S}^2=\id_A$ and both pairs $({\cal A}_L,S)$
and $({\cal A}_L,{\tilde S})$ satisfy (\ref{axi}) and
(\ref{luap}). The map $\gamma_L$ is a composition of a 
coassociative coproduct $\Delta:A\to A\stackrel{\ot}{_{_k}} A$ and
the canonical projection $\p_L: A\stackrel{\ot}{_{_k}} A \to
A\stackrel{\ot}{_{_L}} A$. The compatibility relations
\[ \Delta \ci S =(S\ot S)\ci \Delta^{op} \quad and  \quad
          \Delta \ci {\tilde S} =(S\ot {\tilde S})\ci \Delta^{op} \]
hold true. Then the pair $({\cal A}_L,{\tilde S})$  is a Hopf algebroid
in the sense of Definition \ref{def}.
\ep
\pr We leave to the reader to check that the condition (\ref{luiv})
-- hence all requirements of Proposition \ref{defseq}.{\em ii)} --
hold true.

\section{Integral theory and the dual Hopf algebroid}
\lb{int}
\setc{0}

In this section we generalize the notion of non-degenerate integrals
in (weak) bialgebras to bialgebroids. We examine the consequences of
the existence of a non-degenerate integral in a Hopf algebroid. We do
{\em not} address the question, however, under what conditions on the
Hopf algebroid does the the existence of a non-degenerate integral
follow. That is we do not give a generalization of the Larson-Sweedler
theorem on bialgebroids and neither of the (weaker) Theorem 3.16 in
\cite{BNSz} stating that a weak Hopf algebra possesses a non-degenerate
integral if and only if it is a Frobenius algebra. (About the
implications in one direction see however Theorem 6.3 in \cite{B} and
Theorem \ref{frob} below, respectively.)

{\em Assuming} the existence of a non-degenerate integral in a Hopf
algebroid we show that the underlying bialgebroids are finite. The
duals of finite bialgebroids w.r.t. the base rings were shown to have
bialgebroid structures \cite{KSz} but there is no obvious way how to
transpose the antipode to (either of the four) duals. As the main
result of this section we show that {\em if there exists a
non-degenerate integral} in a Hopf algebroid then the four dual
bialgebroids are all (anti-) isomorphic and they can be made Hopf
algebroids. This dual Hopf algebroid structure is unique up to
isomorphism (in the sense of Definition \ref{hgdeq}).

For the considerations of this section the ``symmetric definition'' of
Hopf algebroids i.e. the characterization in  {\em iii)} of
Proposition \ref{defseq} is the most appropriate. Throughout the
section we use the symmetrized form of the Hopf algebroid introduced
at the end of Subsection \ref{def}.

It is important to emphasize that although the Definitions \ref{intdef}
and \ref{nddef} are formulated in terms of a particular symmetrized
Hopf algebroid ${\cal A}=({\cal A}_L,{\cal A}_R,S)$, actually they
depend only on the Hopf algebroid $({\cal A}_L,S)$. That is to say, if
$\ell$ is a (non-degenerate) left integral in a symmetrized Hopf
algebroid then it is one in any other symmetrized form of the same Hopf
algebroid. Therefore $\ell$ can be called a (non-degenerate) left
integral of the Hopf algebroid. Analogously, although the anti-automorphism
$\xi$ in Lemma \ref{xi} is defined in terms of a  particular
symmetrized Hopf algebroid, it is invariant under the change of the
underlying right bialgebroid.

For a symmetrized Hopf algebroid ${\cal A}=({\cal A}_L,{\cal A}_R,S)$
we use the notations of Section \ref{preli}: the $\atld$ and $\asld$
are the $L$-duals of ${\cal A}_L$, the $\asrd$ and $\atrd$ the
$R$-duals of ${\cal A}_R$. Also for the coproducts of ${\cal A}_L$ and
${\cal A}_R$ we write $\gamma_L(a)=a\di\ot a\dii$ and
$\gamma_R(a)=a\ui\ot a\uii$, respectively. 

\subsection{Non-degenerate integrals}

\bd \lb{intdef} The {\em left integrals} in a left bialgebroid ${\cal
A}_L=(A,L,s_L,t_L,\gamma_L,\pi_L)$ are the invariants of the left
regular $A$ module: 
\[ {\cal I}^L({\cal A})\colon = \{ \ell\in A \vert a\ell=s_L\circ
\pi_L(a)\ell \quad \forall a\in A \} .  \]
The {\em right integrals} in a right bialgebroid ${\cal
A}_R=(A,R,s_R,t_R,\gamma_R,\pi_R)$ are the 
invariants of the right regular $A$ module:
\[ {\cal I}^R({\cal A})\colon = \{ \Upsilon\in A \vert \Upsilon
a=\Upsilon s_R\circ 
\pi_R(a)\quad  \forall a\in A\}. \]
The left/right integrals in a symmetrized Hopf algebroid 
$({\cal A}_L,{\cal A}_R,S)$
are the left/right integrals in ${\cal A}_L/{\cal A}_R$.
\ed  
\bl For a symmetrized Hopf algebroid ${\cal A}$ the following properties of the
element $\ell$ of $A$ are equivalent: 
\[ \begin{array}{rl}
i)&\ell\in  {\cal I}^L({\cal A})\nn
ii)&a\ell = t_L\circ \pi_L(a)\ell {for\ all\ } a\in A\nn
iii)& S(\ell)\in {\cal I}^R({\cal A})\nn
iv)& S^{-1} (\ell)\in {\cal I}^R({\cal A})\nn
v) & S(a) \ell^{(1)} \otimes \ell^{(2)} = \ell^{(1)} \otimes a \ell^{(2)} 
{\ as\  elements\  of\ } A^R \otimes ^RA,  {\ for\ all\ } a\in A.
\nonumber\end{array}\]
\el
\pr Left to the reader.

\bd \lb{nddef}
The left integral $\ell$ in the symmetrized Hopf algebroid ${\cal A}$ is {\em
non-degenerate} if the maps 
\bea \fsr:\asrd\to A &\qquad & \phisr\ \mapsto \  \phisr\ru\ell \ {\rm
and}\nn
\ftr\ :\atrd\to A &\qquad & \phitr\ \, \mapsto \ \phitr\rd \ell
\eea 
are bijective. The right integral $\err$  in the symmetrized Hopf algebroid ${\cal
A}$ is
{\em non-degenerate} if $S(\err)$ is a non-degenerate left integral
i.e. if the maps 
\bea \err_L:\atld\to A &\qquad & \phitl\ \mapsto \ \err\lu \phitl \
\quad {\rm and}\nn
{_L\err}:\asld\to A &\qquad &\phisl\ \,\mapsto \ \err\ld \phisl
\eea
are bijective.
\ed
\br \lb{cop}
If $\ell$ is a non-degenerate left integral in the symmetrized Hopf algebroid
${\cal A}$ then so is 
in ${\cal A}_{cop}$, and when replacing ${\cal A}$
with ${\cal A}_{cop}$  the roles of $\fsr$ and $\ftr$
become interchanged. Hence 
any statement proven in a symmetrized Hopf algebroid possessing a non-degenerate
left integral on $\fsr$ implies that the co-opposite statement holds
true on $\ftr$. 
\er 
\bt \lb{frob} 
Let ${\cal A}$ be a  symmetrized Hopf algebroid possessing a
non-degenerate left integral. Then the ring extensions 
$s_R:R\to A$, $t_R:R^{op}\to A$, $s_L:L\to A$ and $t_L:L^{op}\to A$
are all Frobenius extensions.
\et 
\pr Let $\ell$ be a non-degenerate left integral in ${\cal A}$. With
its help we construct the Frobenius system for the extension $s_R:R\to
A$. It consists of a Frobenius map 
\be \lsr\colon = \fsr^{-1}(1_A) :{_R A^R} \to R \ee
and a quasi-basis (in the sense of (\ref{quasib})) for it:
\[ \ell\ui\ot S(\ell\uii)\in \asr \ot {_R A}.\]
As a matter of fact the $\lsr$ is a right $R$-module map $\asr\to R$
by construction. We claim that it is also a left $R$-module map ${_R
A}\to R$. Since
\be  (\lsr\lu S(a))\ru \ell= 
\ell\uii t_R\ci \lsr \left(S(a)\ell\ui\right)=
a(\lsr\ru \ell )=a\lb{fsrinv}\ee
the inverse $\fsr^{-1}$ maps $a\in A$ to $\lsr \lu S(a)$. This implies
that $\lsr \lu s_R(r)=\fsr^{-1}\ci t_R(r)$. Now for a given element
$r\in R$ the map $\chi(r)^*:\asr \to R$, $a\ \mapsto \ r\lsr(a)$ is also
equal to $\fsr^{-1}\ci t_R(r)$:
\[ \chi(r)^*\ru \ell= \ell\uii t_R(r\lsr(\ell\ui))=(\lsr\ru
\ell)t_R(r)=t_R(r).\]
Applying the two equal maps $\lsr\lu s_R(r)$ and $\chi(r)^*$ to an
element $a\in A$ we obtain
\be \lsr(s_R(r)a)=r\lsr(a).\lb{srl} \ee
This proves that $\lsr$ is an $R$-$R$ bimodule map ${_R A^R}\to A$. Also 
\bea s_R\ci \lsr (a\ell\ui) S(\ell\uii)&=&
   S( \lsr\ru \ell)a=a \quad {\rm and}\lb{qb1}\\
\ell\ui s_R\ci \lsr( S(\ell\uii)a)&=& a  \ell\ui s_R\ci \lsr\ci
S(\ell\uii).\eea
Now we claim that $\ell\ui s_R\ci \lsr\ci S(\ell\uii)\equiv \lsr\ci
S\rd \ell$ is equal to $1_A$, which proves the claim. (Recall that by
(\ref{srl}) $\lsr\ci S$  is an element of $\atrd$.)
Since  $\fsr^{-1}(a)=\lsr \lu S(a)$, 
\be  \phisr(a)=[\lsr \lu S(\phisr\ru \ell)](a)=
\lsr\left(s_R\ci \phisr(\ell\ui)S(\ell\uii)a\right)\lb{rn}\ee
for all $\phisr\in \asrd$ and $a\in A$.

By the bijectivity of $\ftr$ we can introduce the element $\ltr\colon
=\ftr^{-1} (1_A)\in \atrd$. Analogously to (\ref{srl}) and
(\ref{fsrinv}) we have
\bea \ltr(t_R(r)a)&=&\ltr(a)r \quad {\rm and}\lb{trl}\\
     \ftr^{-1}(a)&=&\ltr\ld S^{-1}(a).\eea
Using the fact that both $\ftr$ and $S^{-1}$ 
are bijective so is the map $A\to \atrd$, $a\ \mapsto\  \ltr\ld
a$. Using the identities  (\ref{trl}), (\ref{srl}) and (\ref{rn}) compute
\bea 
\left( \ltr \ld S^{-1}(\lsr\ci S \rd \ell)\right)(a)&=&
\ltr\left( t_R\ci \lsr \ci S(\ell\uii)S^{-1}(\ell\ui) a\right)=\nn
&=&\lsr \left(s_R\ci \ltr \ci S^{-1}(\ell\ui)S(\ell\uii) S(a)\right)=
\ltr(a)\nonumber\eea
for all $a\in A$. This is equivalent
to $\lsr \ci S\rd \ell=1_A$ that is $\ltr=\lsr\ci S$ 
proving that $(\lsr, \ell\ui\ot
S(\ell\uii))$ is a Frobenius system for the extension $s_R:R\to A$.

By repeating the same proof in ${\cal A}_{cop}$ we obtain the Frobenius
system $(\ltr, \ell\uii\ot S^{-1}(\ell\ui))$ for the extension
$t_R:R\to A$. 

It is straightforward to check that $(\mu^{-1}\ci \lsr, \ell\ui \ot
S(\ell\uii))$ is a Frobenius system for the extension $t_L:L\to A$ and
$(\nu^{-1}\ci \ltr, \ell\uii\ot S^{-1}(\ell\ui))$ is a Frobenius
system for the extension $s_L:L\to A$.
\qed

\smallskip

From now on let $\ell$ be a non-degenerate left integral in the
symmetrized Hopf algebroid ${\cal A}$, set
$\lsr=\fsr^{-1}(1_A)$ and $\ltr=\ftr^{-1} (1_A)$.

Theorem \ref{frob} implies that for a symmetrized Hopf algebroid ${\cal A}$ 
possessing a non-degenerate integral the modules $\asr$, $\atr$, $\atl$  and
$\asl$ are finitely 
generated projective. Hence by the result of \cite{KSz}, their duals
$\asrd$ and $\atrd$ carry left bialgebroid 
structures over the base $R$,  and $\atld$  and
$\asld$ carry right bialgebroid structures over the base $L$:
\be\begin{array}{ll}
     s^*_L(r)(a)= r \pi_R(a) \qquad &^*s_L(r)(a)= \pi_R(t_R(r)a) \nn
     t^*_L(r)(a)=\pi_R(s_R(r)a) \qquad &^*t_L(r)(a)=\pi_R(a) r \nn
     \gamma^*_L(\phisr)= \phisr\lu \ell\ui \ot \fsr^{-1}(\ell\uii)\qquad
     &^*\gamma_L (\phitr)= \ftr^{-1}(\ell\ui) \ot \phitr\ld \ell\uii\nn
     \pi^*_L(\phisr)= \phisr(1_A) \qquad &^*\pi_L(\phitr)= \phitr(1_A)\nn
\quad\nn
     s_{*R}(l)(a)=\pi_L(a s_L(l)) \qquad &_*s_R(l)(a)= \pi_L(a)l \nn
     t_{*R}(l)(a)=l \pi_L(a)  \qquad &_*t_R(l)(a)=\pi_L(a t_L(l))\nn
     \gamma_{*R}(\phitl)=\ell\di \ru \phitl \ot \ftl^{-1}(\ell\dii) 
     \qquad
     &_*\gamma_R(\phisl)=\fsl^{-1}(\ell\di)\ot \ell\dii\rd \phisl\nn
     \pi_{*R}(\phitl)= \phitl(1_A) \qquad  & _*\pi_R(\phisl)=\phisl(1_A)
\end{array}\lb{dbgd}\ee

\bl \lb{lac}
Let $\ell$ be a non-degenerate left integral in the symmetrized Hopf 
algebroid ${\cal A}$. Then for $\lsr=\fsr^{-1}(1_A)$, $\ltr=\ftr^{-1} 
(1_A)$ and any element $a\in A$ the identities 
\bea \lsr \rightharpoonup a &=& s_R \circ \lsr (a) \lb{lsac}\\
 \ltr \rightharpoondown a &=& t_R \circ \ltr (a)  \lb{ltac}
\eea
hold true.
\el
\pr One checks that 
\[ \phisr\lsr=\fsr^{-1} ( \phisr \rightharpoonup 1_A) = s^*_L \ci
\phisr(1_A) \lsr \] 
for all $\phisr\in \asrd$. This implies that $\phisr (\lsr \ru
a)=\phisr (s_R \ci \lsr (a) )$ for all $\phisr\in \asrd$. Since $\asr$ is
finitely generated projective by Theorem \ref{frob} this proves
(\ref{lsac}). The identity (\ref{ltac}) follows by Remark \ref{cop}.
\qed

\medskip

The left integrals
in a Hopf algebroid were defined in Definition \ref{intdef} as the left
integrals in the underlying left bialgebroid. The non-degeneracy of
the left integral was defined in Definition
\ref{nddef} using however the underlying right bialgebroid as well,
that is it relies to  the whole of the Hopf algebroid
structure. Therefore it is not obvious whether 
the non-strict isomorphisms of Hopf algebroids preserve non-degenerate
integrals. In  the rest of this subsection we prove that this is the
case: 

\bp \lb{ndtw}
Let both $({\cal A}_L,S)$ and 
$({\cal A}_L,S^{\prime})$ be
Hopf algebroids.
Then their  non-degenerate left integrals  coincide.
\ep
\pr A left integral $\ell$ in $({\cal A}_L,S)$ is a left integral in
$({\cal A}_L,S^{\prime})$ by definition.  

Let ${\cal A}_R=(A,R,s_R,t_R,\gamma_R,\pi_R)$ and ${\cal
A}_R^{\prime}=(A,R^{\prime},s_R^{\prime},t_R^{\prime},\gamma_R^{\prime},
\pi_R^{\prime})$ be the right bialgebroids underlying the Hopf
algebroids $({\cal A}_L,S)$ and $({\cal A}_L,S^{\prime})$, respectively.
It follows from the uniqueness of the maps $\alpha^{-1}$ and
$\beta^{-1}$ in (\ref{schinv}) that the coproducts
$\gamma_R(a)=a\ui\ot a\uii$ of ${\cal A}_R$ and $\gamma^{\prime}_R(a)=
a^{\{1\}}\ot a^{\{2\}} $ of ${\cal A}^{\prime}_R$ are related as
\be a\ui \ot S^{\prime -1}\ci S(a\uii)= a^{\{1\}}\ot
a^{\{2\}} = S^{\prime}\ci S^{-1}(a\ui)\ot a\uii.
\lb{twcp}\ee
With the help of the maps $\mu=\pi_R\ci t_L$,
$\mu^{\prime}=\pi^{\prime}_R\ci t_L$, $\nu=\pi_R\ci s_L$ and
$\nu^{\prime}=\pi^{\prime}_R\ci s_L$  we can introduce the isomorphisms
of additive groups
\bea \asrd&\to\  {\cal A}^{\prime *}\qquad \phisr&\mapsto\  \mu^{\prime}\ci
\mu^{-1}\ci \phisr\nn
\atrd&\to\  {^*{\cal A}^{\prime}}\qquad \phitr &\mapsto\  \nu^{\prime}\ci
\nu^{-1}\ci \phitr.\nonumber\eea
Then the canonical actions (\ref{srac}) of $ {\cal A}^{\prime *}$ and
$\asrd$ and of ${^*{\cal A}^{\prime}}$ and $\atrd$ on $A$ are related as
\bea  \mu^{\prime}\ci \mu^{-1}\ci \phisr \stackrel{\prime}{\ru} a
&=&S^{\prime -1}\ci S(\phisr \ru a)\nn 
\nu^{\prime}\ci \nu^{-1}\ci \phitr\stackrel{\prime}{\rd} a &=& S^{\prime}\ci
S^{-1}(\phitr\rd a)\nonumber\eea
what implies the non-degeneracy of the left integral $\ell$ in $({\cal
A}_L,S^{\prime})$  provided it is non-degenerate in $({\cal A}_L,S)$.
\hspace{1cm}\qed


\subsection{Two sided non-degenerate integrals}

The Proposition \ref{ndtw} above proves that the structure of the
non-degenerate left integrals is the same within an isomorphism class
of Hopf algebroids. In this subsection we prove that for a
non-degenerate left integral $\ell$ in the Hopf algebroid $({\cal
A}_L,S)$ there exists a distinguished representative $({\cal
A}_L,S^{\prime}_{\ell}) $ in the
isomorphism class of $({\cal A}_L,S)$ with the property that $\ell$ is
not only a non-degenerate left integral in $({\cal
A}_L,S^{\prime}_{\ell})$ but also a non-degenerate right integral. 

The Hopf algebroids with two sided non-degenerate integral are of
particular interest. Both the Hopf algebroid structure constructed on
the dual of a Hopf algebroid in Subsection \ref{dual} and the one
associated to a depth 2 Frobenius extension in Section \ref{ringex}
belong to this class.

\bl \lb{lprime}
Let $\ell$ be a non-degenerate left integral in a symmetrized Hopf algebroid 
${\cal A}$. Set $\lsr\colon =\fsr^{-1}(1_A)$ and $\ltr\colon
=\ftr^{-1}(1_A)$. Then any (not necessarily non-degenerate) left
integral $\ell^{\prime} \in {\cal I}^L({\cal A})$ satisfies
\[ \ell s_R\ci \lsr (\ell^{\prime}) =\ell^{\prime} = 
    \ell t_R\ci \ltr (\ell^{\prime}).\]
\el
\pr Observe that for $\phitr\in\atrd$ and $\ell^{\prime}\in {\cal
I}^L({\cal A})$ we have 
$\phitr\ld  S^{-1} (\ell^{\prime}) = {^*\!t_L} \ci\phitr \ci S^{-1}
(\ell^{\prime})$ hence 
\[ \ell^{\prime} = (\ltr \ld S^{-1} (\ell^{\prime}))\rd \ell =
{^*\!t_L} \ci\ltr \ci S^{-1} (\ell^{\prime})\rd \ell = \ell s_R\ci \lsr
(\ell^{\prime} ).
\]
The identity $\ell^{\prime} = \ell t_R\ci \ltr (\ell^{\prime})$
follows by Remark \ref{cop}.
\qed

\bl \lb{xi} 
Let $\ell$ be a non-degenerate left integral in a symmetrized Hopf algebroid 
${\cal A}$. Set $\lsr\colon =\fsr^{-1}(1_A)$.
Then the map $\xi: A\to A$; $a\mapsto S((\lsr \lu \ell)\ru a)$ is a ring
anti-automorphism. 
\el
\pr  Using Lemma \ref{lprime} one checks that
\bea [(\lsr \lu \ell)\ru a]&&\!\!\!\!\!\!\!\!\!\!\!\!\![(\lsr \lu \ell)\ru b]=
a\uii t_R\ci \lsr(\ell a\ui) b\uii t_R\ci \lsr(\ell b\ui)=\nn
&=&a\uii b\uii t_R\ci \lsr \left(\ell s_R\ci\lsr(\ell
a\ui)b\ui\right)=
a\uii b\uii t_R\ci \lsr \left(\ell a\ui b\ui\right)=\nn
&=&(\lsr \lu \ell)\ru ab
\nonumber\eea
for $a,b\in A$, hence  the map $\xi$ is anti-multiplicative. By analogous
calculations the reader may check that it is bijective with inverse 
$\xi^{-1}(a)=S^{-1} ((\ltr \leftharpoondown \ell)\rightharpoondown
a)$, where  $\ltr\colon =\ftr^{-1}(1_A)=\lsr\ci S$.
\hspace{1cm}\qed 

\bp \lb{flbij} 
Let $\ell$ be a non-degenerate left integral in a symmetrized Hopf algebroid 
${\cal A}$. Then the maps 
\[\begin{array}{cccc}
 \ftl: \atld \to A  &\qquad \phitl&\mapsto& \quad\ell\leftharpoonup
\phitl \nn
\fsl: \asld \to A  &\qquad   \phisl&\mapsto& \quad \ell\leftharpoondown
\phisl
\end{array}\]
are bijective.
\ep
\pr We claim that with the help of the ring isomorphism $\nu$
(introduced in Proposition \ref{sbgdi}) we have $\ell\ld \phisl =
\xi^{-1}\circ \fsr (\nu\ci \phisl \ci S^{-1})$, which implies the
bijectivity of $\fsl$. As a matter of fact 
\bea \ell& \ld& _*[\nu^{-1} \ci \fsr^{-1}(a)\ci S] =
t_L\ci \nu^{-1} \ci (\lsr \lu S(a))\ci S (\ell\dii)\ell\di =
S^{-1}\ci t_R \ci \ltr (\ell\dii a) \ell\di=\nn
&=& S^{-1}(\ltr \rd \ell\dii a) \ell\di=
S^{-1}\left( S(\ell\di) {\ell\dii}\ui a\ui s_R \ci \ltr 
({\ell\dii}\uii a\uii)\right)= \nn
&=&  S^{-1}\left( s_R\ci\pi_R (\ell\ui) a\ui s_R \ci \ltr (\ell\uii
a\uii)\right)=
  S^{-1}\left( a\ui s_R \ci \ltr (\ell a\uii)\right)=\nn
&=&  S^{-1}((\ltr \ld \ell)\rd a)=\xi^{-1}(a).
\lb{fsl}\eea 
Similarly, by the application of (\ref{fsl}) to ${\cal A}_{cop}$
\be \ell\lu \phitl =\xi\ci \ftr (\mu\ci \phitl \ci S), \lb{ftl}\ee
 hence $\ftl$ is also bijective.
\hspace{1cm} \qed

Using (\ref{ftl}) we have an equivalent
form of the anti-automorphism $\xi$ introduced in Lemma \ref{xi}:
\be \xi(a)= 
\ell\lu (a\ru {\ell_L}^{-1}(1_A)). \lb{xiarr}\ee
\bl \lb{arrrel}
Let $\ell$ be a non-degenerate left integral in a symmetrized Hopf algebroid 
${\cal A}$.
Then for all elements $a,b\in A$ we have the identities
\bea \fsr^{-1}(b)\ru a &=& \ftr^{-1}(a)\rd b \lb{srtr}\\
     a\lu\ftl^{-1}(b) &=& b\ld \fsl^{-1}(a) \lb{sltl}\\
     \fsr^{-1}(b)\ru a &=& a\lu \ftl^{-1}(b)\lb{srtl}\\
     \ftr^{-1}(b)\rd a &=& a\ld \fsl^{-1}(b).\lb{sltr}
\eea
\el
\pr We illustrate the proof on (\ref{srtr}). Use Lemma \ref{lac} to
see that
\[ \ftr^{-1}(a) \rd b =
b\ui \left(\lsr\ru S(b\uii)a\right)=
s_L\ci\pi_L(b\di)a\uii t_R\ci \lsr \left(S(b\dii)a\ui\right)=
\fsr^{-1}(b) \ru a\]
where $\lsr=\fsr^{-1}(1_A)$. The rest of the proof is analogous. 
\hspace{1cm}\qed

\bl \lb{kappa}
Let $\ell$ be a non-degenerate left integral in a symmetrized Hopf algebroid 
${\cal A}$. Set  $\lsr=\fsr^{-1}(1_A)$. Then the map $\kappa:R\to R$,
$r\mapsto  \lsr(\ell t_R(r))$ is a 
ring automorphism.
\el 
\pr It follows from Lemma \ref{lprime} that  $\kappa$ is
multiplicative: for $r,r^{\prime}\in R$
\[ \kappa(r) \kappa (r^{\prime})  =
\lsr(\ell t_R(r)) \lsr(\ell t_R(r^{\prime}))=
\lsr (\ell s_R\ci \lsr(\ell t_R(r^{\prime})) t_R(r))=
\lsr (\ell t_R(r^{\prime}) t_R(r) )= \kappa(r r^{\prime}).
\]
In order to show that $\kappa$ is bijective we construct the inverse
$\kappa^{-1} : r\mapsto \ltr(\ell s_R(r))$ where
$\ltr=\ftr^{-1}(1_A)=\lsr\ci S$. 
\hspace{1cm}\qed

\bp \lb{aprl} Let $\ell$  be a non-degenerate left integral in the
Hopf algebroid $({\cal A}_L,S)$. Then there 
exists a unique Hopf algebroid $({\cal A}_L,S^{\prime}_{\ell})$
such that $\ell$ is a two sided
non-degenerate integral in $({\cal A}_L,S^{\prime}_{\ell})$.
\ep

\pr {\em uniqueness:} Suppose that $({\cal
A}_L,S^{\prime}_{\ell})$ is a Hopf algebroid of the required
kind. Denote the underlying right bialgebroid by ${\cal A}^{\prime}_R=
(A,R^{\prime},s_R^{\prime}, t_R^{\prime}, \gamma_R^{\prime},
\pi_R^{\prime},)$. Define $\lambda^{\prime *}\in {\cal A}^{\prime 
*}$  with the property that $\lambda^{\prime
*}\stackrel{\prime}{\ru}\ell =1_A$ (where $\stackrel{\prime}{\ru}$
denotes the canonical action (\ref{srac}) of ${\cal A}^{\prime *}$ on
$A$). Introducing the notation $\gamma_R^{\prime}(a)=a^{\{1\}}\ot
a^{\{2\}}$  one checks that
\[ S^{\prime-1}(\ell)=(\lambda^{\prime *}\stackrel{\prime}{\lu} \ell)
\stackrel{\prime}{\ru} \ell=
\ell^{\{2\}}t^{\prime}_R\ci \lambda^{\prime *}(\ell
\ell^{\{1\}}) =
\ell^{\{2\}}t^{\prime}_R\ci \lambda^{\prime *}\left(\ell
s^{\prime}_R\ci \pi^{\prime}_R(\ell^{\{1\}})\right) =\nn
\ell t^{\prime}_R\ci \lambda^{\prime *}(\ell)=\ell. \]
With the help  of the
element $\pi_L\ci S^{-1}\ci S^{\prime}\in \atld$ 
we have 
\[ S(a\lu \pi_L\ci S^{-1}\ci S^{\prime})=
S\ci S^{\prime-1}\left(S^{\prime}(a)^{\{1\}}\right)s_R\ci
\pi_R\left(S^{\prime}(a)^{\{2\}}\right)=
S^{\prime}(a)\]
for all $a\in A$ where in the last step the relation (\ref{twcp}) has
been used. Then the condition $S^{\prime}(\ell)=\ell$ is equivalent to
\[ S^{\prime}(a)=S(a\lu \ell_L^{-1}\ci S^{-1}(\ell)).\]
This proves the uniqueness of $S^{\prime}$.

{\em existence:} Let $\xi$ be the anti-automorphism of $A$ introduced
in Lemma \ref{xi}.
We claim that $({\cal A}_L,\xi)$ is a Hopf algebroid
of the required kind. Introduce the  right
bialgebroid ${\cal A}^{\prime}_R$ on the total ring $A$ over the base
$R$ with  structural maps  
\[ s^{\prime}_R=s_R \qquad   t^{\prime}_R=\xi^{-1}\ci s_R \qquad
\gamma^{\prime}_R=\xi^{-1}_{A\ot_L A}\ci S_{A\ot_R A}\ci \gamma_R\ci
S^{-1}\ci \xi 
\qquad 
\pi^{\prime}_R=\pi_R\ci S^{-1}\ci \xi \] 
where ${\cal A}_R=(A,R,s_R,t_R,\gamma_R,\pi_R)$ is the right
bialgebroid underlying $({\cal A}_L,S)$. First we check that
the triple $({\cal A}_L,{\cal A}^{\prime}_R,\xi)$ satisfies
Proposition \ref{defseq} {\em iii)}. Since
\[ \xi^{-1}\ci S\ci s_R=\xi^{-1}\ci t_R\ci \theta_R=S^{-1}\ci t_R\ci
\theta_R=s_R \quad {\rm and} \quad 
\xi^{-1}\ci S\ci t_R=\xi^{-1}\ci s_R \]
the ${\cal A}^{\prime}_R$ is a right bialgebroid isomorphic to ${\cal
A}_R$ via  the isomorphism $(\xi^{-1}\ci S,\id_R)$.

The requirement $s^{\prime}_R(R)\equiv s_R(R)=t_L(L)$ obviously holds
true. Since 
\[ t^{\prime}_R(r)=\xi^{-1}\ci s_R(r)= S^{-1}\ci s_R\ci\ltr(\ell s_R(r))=
t_R\ci \kappa^{-1}(r)\]
also $t^{\prime}_R(R)\equiv t_R(R)=s_L(L)$.
Since
\bea \gamma^{\prime}_R(a)&\equiv& a^{\{1\}}\ot a^{\{2\}}=
\xi^{-1}_{A\ot_L A}\ci\gamma_L\ci \xi(a)=
a\ui \ot \xi^{-1}\ci S(a\uii)=\nn
&=&a\ui \ot S^{-1}\left(\ftr^{-1}\ci S(\ell)\rd S(a\uii)\right)=
a\ui \ot S^{-1}\left(S(a\uii)\ld {_L\ell}^{-1}\ci S(\ell)\right)=\nn
&=&a\ui s_R\ci \kappa\ci \mu\ci {_L\ell}^{-1}\ci S(\ell)\ci S(a\uii)\ot a^{(3)}
\nonumber\eea
we have
\bea (\gamma_L\ot \id_A)\ci \gamma^{\prime}_R(a)&=&
{a\ui}\di \ot {a\ui}\dii\ot \xi^{-1}\ci S(a\uii)=
a\di\ot {a\dii}\ui\ot \xi^{-1}\ci S({a\dii}\uii)=\nn
&=&(\id_A\ot \gamma^{\prime}_R)\ci \gamma_L(a)\nn
(\id_A \ot \gamma_L)\ci \gamma^{\prime}_R(a) &=&
a\ui s_R\ci \kappa\ci \mu \ci {_L\ell}^{-1}\ci S(\ell)\ci S(a\uii)\ot {a^{(3)}}\di 
\ot a^{(3)}\dii=\nn
&=&{a\di}\ui s_R\ci \kappa\ci \mu\ci {_L\ell}^{-1}\ci S(\ell)\ci S({a\di}\uii)\ot
{a\di}^{(3)}\ot a\dii=\nn
&=&(\gamma^{\prime}_R\ot \id_A)\ci \gamma_L(a).
\nonumber\eea
By Lemma \ref{xi} the $\xi$ is an anti-automorphism of the ring $A$. The identity 
$\xi\ci t^{\prime}_R=s^{\prime}_R$ is obvious and also
\[ \xi\ci t_L=\xi\ci s_R\ci\mu=S\ci s_R\ci \mu=s_L.\]
Finally,
\bea \xi(a\di)a\dii&=&
s_R\ci \pi_R\left((\lsr\lu \ell)\ru a\right)=
s_R\ci \pi_R\ci S^{-1}\ci \xi(a)=
s^{\prime}_R\ci \pi^{\prime}_R (a)\nn
a^{\{1\}}\xi(a^{\{2\}})&=&
a\ui \xi\ci \xi^{-1}\ci S(a\uii)=s_L\ci \pi_L(a).
\nonumber\eea
This proves that ${\cal A}^{\prime}_{\ell}=({\cal A}_L,{\cal
A}^{\prime}_R,\xi)$  satisfies Proposition \ref{defseq} {\em
iii)} hence $({\cal A}_L,\xi)$ is a Hopf algebroid. Since   
\[\xi(\ell)=S\left((\lsr\lu \ell)\ru \ell\right)=S\ci S^{-1}(\ell)=\ell\]
the $\ell$ is a two sided non-degenerate integral in ${\cal
A}^{\prime}_{\ell}$. 
\hspace{1cm}\qed


\subsection{Duality}
\lb{dual}

It follows from Theorem \ref{frob} that for a symmetrized Hopf
algebroid ${\cal A}$ possessing 
a non-degenerate left integral $\ell$ the dual rings (with respect to
the base 
ring) carry bialgebroid structures. These bialgebroids (\ref{dbgd}) are
independent 
of the particular choice of the non-degenerate integral. In this
subsection we  
analyze these bialgebroids. We show that the four bialgebroids
(\ref{dbgd}) are 
all (anti-) isomorphic and can be equipped with an $\ell$-dependent
Hopf algebroid  
structure. 
Because of the  $\ell$-dependence of this Hopf algebroid structure the
duality of Hopf algebroids is sensibly defined on the isomorphism
classes of Hopf algebroids. 

\bl Let $\ell$ be a non-degenerate left integral in the symmetrized Hopf algebroid
${\cal A}$. Then with the help of the anti-automorphism $\xi$ of Lemma
\ref{xi} we have the equalities 
\be  \xi^{-1}(\ell\uii)\ot S^{-1}(\ell\ui) = \ell\di\ot \ell\dii =
S(\ell^{\uii})\ot \xi (\ell\ui) \lb{lud}\ee
in $\atl\ot \asl$.
\el
\pr  The element $\xi^{-1}(\ell\uii)\ot S^{-1}(\ell\ui)$ is in $\atl\ot
\asl$  since 
$S^{-1}\ci s_R =t_R= s_L\ci \nu^{-1}$ and $\xi^{-1}\ci t_R =S^{-1}\ci t_R =
t_L\ci \nu^{-1}$. Using (\ref{fsl}), in $\atl\ot \asl$ we have 
\[ \xi^{-1}(\ell\uii)\ot S^{-1}(\ell\ui) =
\ell\di \ot t_R \ci \ltr (\ell\dii \ell\uii) S^{-1}(\ell\ui)=
\ell\di \ot \ell\dii.\]
The other equality follows by repeating the proof in ${\cal A}_{cop}$.
\hspace{1cm}\qed

\bc For a  non-degenerate left integral $\ell$ in the symmetrized Hopf algebroid
${\cal A}$ the maps $\ftl$ and $\fsl$ satisfy the identities
\bea \ftl(a\ru \phitl) &=& \ftl(\phitl)
\xi(a) \lb{ftlint}\\
     \fsl(a\rd \phisl) &=& \fsl(\phisl) \xi^{-1}(a)
\lb{fslint}\eea
where $\xi$ is the anti-automorphism of $A$ introduced in Lemma
\ref{xi}. 
\ec

\bt \lb{bgdiso}
Let $\ell$ be a non-degenerate left integral in the symmetrized Hopf algebroid
${\cal A}=({\cal A}_L,{\cal A}_R,S)$. Then
the left bialgebroids $\asrd_L$, $\atrd_L$, $(\atld_R)^{op}_{cop}$
and $(_*{\cal A}_R)^{op}_{cop}$ in (\ref{dbgd}) are isomorphic  via the
isomorphisms 
\begin{center}
\begin{picture}(252,80)
\put(67,0){$\asrd_L$}
\put(214,0){$\atrd_L$}
\put(59,68){$(\atld_R)^{op}_{cop}$}
\put(207,68){$(\asld_R)^{op}_{cop}$}
\put(90,2){\vector(1,0){120}}
\put(100,70){\vector(1,0){105}}
\put(70,65){\vector(0,-1){55}}
\put(217,65){\vector(0,-1){55}}
\put(0,32){$(\fsr^{-1}\ci\ftl,\nu)$}
\put(222,32){$(\ftr^{-1}\ci \fsl,\mu)$}
\put(100,8){$(\ftr^{-1}\ci \xi^{-1}\ci \fsr,\theta_R^{-1})$}
\put(100,76){$(\fsl^{-1}\ci\xi^{-1} \ci \ftl,\id_R)$}
\end{picture}
\end{center}
where $\xi$ is the anti-automorphism of $A$ introduced in Lemma
\ref{xi} and the maps $\mu,\nu$ and $\theta_R$ are the ring isomorphisms
introduced in Proposition \ref{sbgdi}.
\et

\pr  By Proposition \ref{sbgdi} the map $\nu$ is a ring isomorphism
$L^{op}\to R$. By Proposition \ref{flbij} $\fsr^{-1}\ci \ftl$ is
bijective. Its anti-multiplicativity follows from (\ref{srtl}).
The comultiplicativity follows by the successive use of the identity 
$\fsr(\phisr\lu a)=S^{-1}(a)\fsr(\phisr)$, the integral property of $\ell$,
(\ref{lud})  and (\ref{ftlint}):
\bea \gamma^{*}_L&\ci& \fsr^{-1}\ci \ftl (\phitl) =
\fsr^{-1}\ci \ftl (\phitl)\lu \ell\ui \ot
\fsr^{-1}(\ell\uii) =\nn
&=&\fsr^{-1} \left(S^{-1}(\ell\ui)  \ftl (\phitl)\right) \ot
\fsr^{-1}(\ell\uii) =
 \fsr^{-1} \left(S^{-1}(\ell\ui)\right) \ot
\fsr^{-1}(\ftl (\phitl)\ell\uii) =\nn
&=&\fsr^{-1} (\ell\dii) \ot
\fsr^{-1}\left(\ftl (\phitl)\xi(\ell\di)\right) =
\fsr^{-1}\ci \ftl (\ftl^{-1}(\ell\dii)) \ot
\fsr^{-1}\ci \ftl (\ell\di\ru\phitl)=\nn 
&=& (\fsr^{-1}\ci \ftl \ot \fsr^{-1}\ci \ftl)\ci
{\gamma_{*R}}^{op} (\phitl).
\nonumber\eea
One checks also
\bea \left( \fsr^{-1}\ci \ftl \ci s_{*R}(l) \right)(a) &=& 
\ltr \left( S^{-1}(a) s_L\ci \pi_L[\ell\di s_L(l)]\ell\dii \right)=
\nu(l) \pi_R(a) =
\left( s^*_L\ci \nu(l) \right)(a)\nn
 \left( \fsr^{-1}\ci \ftl \ci t_{*R}(l) \right)(a) &=& 
\ltr \left( S^{-1}(a) s_L[l\pi_L(\ell\di)]\ell\dii \right)=
 \pi_R\left( s_R\ci\nu(l)a \right)= \left(t^*_L\ci
\nu(l)\right)(a)\nn
\pi^*_L \ci \fsr^{-1} \ci \ftl (\phitl) &=&
\nu\ci\phitl\left( t_L\ci \nu^{-1} \ci \ltr(\ell\dii) \ell\di \right)=\nu\ci 
\phitl(1_A)= \nu \ci \pi_{*R}(\phitl).
\nonumber \eea
This proves that $( \fsr^{-1} \ci \ftl,\nu)$ is a
bialgebroid isomorphism $(\atld_R)^{op}_{cop} \to \asrd_L$.
By Remark \ref{cop} ($\ftr^{-1}\ci \fsl,\mu)$ is  a
bialgebroid isomorphism  $(\asld_R)^{op}_{cop}\to \atrd_L$.

By (\ref{ftl}) $\ftr^{-1}\ci \xi^{-1}\ci \fsr =\mu
\ci\ftl^{-1}\ci\fsr(\phisr)\ci S$ hence we have to prove that
$(\mu\ci - \ci S,\mu)$ is a 
bialgebroid isomorphism $({\cal A}_{*R})^{op}_{cop}\to \atrd_L$. 

The map $\mu$ is a ring isomorphism $L^{op}\to R$ by Proposition
\ref{sbgdi}. 
The map $\phitl\mapsto \mu\ci \phitl \ci S$ is bijective.
Its anti-multiplicativity is obvious.
The anti-comultiplicativity follows from (\ref{ftl}) and (\ref{lud}):
\bea ^*\gamma_L (\mu\ci \phitl\ci S) &=&
\ftr^{-1}(\ell\ui) \ot \mu \ci \phitl \ci S \ld \ell\uii= \nn 
&=&\mu \ci \ftl^{-1}\ci\xi(\ell\ui) \ci S \ot 
\mu \ci (S(\ell\uii)\ru \phitl)\ci S  =
\mu\ci \phitl\uii \ci S \ot \mu\ci \phitl\ui \ci S\ .
\nonumber\eea
Finally, by Proposition \ref{sbgdi}
\bea \left(\mu \ci s_{*R} (l) \ci S \right) (a) &=&
\mu \ci \pi_L \left( S(a) s_L(l) \right)=
\pi_R \left( s_R\ci \mu (l) a\right) =
\left({^*\!s_L}\ci \mu(l)\right) (a)\nn
 \left(\mu \ci t_{*R} (l) \ci S \right) (a) &=&
\mu\left(l \pi_L\ci S(a)\right)= \pi_R(a)\mu(l)=
\left({^*\!t_L}\ci \mu(l)\right)(a)\nn
{^*\! \pi_L}\left(\mu\ci \phitl \ci S\right) &=& 
\mu\ci \phitl\ci S(1_A)=
\mu\ci \pi_{*R}(\phitl).
\nonumber \eea
This proves the theorem.
\qed

\bt \lb{dualhgd}
Let $\ell$ be a non degenerate  
left integral in the symmetrized Hopf algebroid ${\cal A}$. Then the left
bialgebroid ${\cal 
A}_{*L}^{\ell}=(\atld,R,s_{*L}, 
t_{*L},\gamma_{*L},\pi_{*L})$ where
\[\begin{array}{ll} 
s_{*L}(r)(a)=\mu^{-1}(r)\pi_L(a)\qquad&
t_{*L}(r)(a)=\pi_L(a t_R\ci \kappa^{-1}(r))\\ 
\gamma_{*L}(\phitl)=\xi^{-2}(\ell\dii)\ru \phitl \ot
{\ell_L}^{-1}(\ell\di)\qquad& 
\pi_{*L}(\phitl)=\lsr(\ell\lu\phitl),
\end{array}\]
the right bialgebroid  ${\cal A}_{*R}$ in (\ref{dbgd}) and the antipode
$S_*^{\ell}:= {\ell_L}^{-1}\ci \xi \ci {\ell_L}$ form a symmetrized Hopf algebroid
denoted by $\atld^{\ell}$.
\et

\pr We show that the triple $\atld^{\ell}\colon =({\cal
A}_{*L}^{\ell},$ ${\cal A}_{*R},S_*^{\ell})$ satisfies the {\em iii)} of
Proposition \ref{defseq}. 

The ${\cal A}_{*L}^{\ell}$ is a left bialgebroid isomorphic to 
$(\asld_R)^{op}_{cop}$ via the isomorphism $({_L\ell}^{-1}\ci
\ell_L,\mu)$. Also 
\[ s_{*L}(R)=t_{*R}(L) \qquad {\rm and}\qquad
   t_{*L}(R)=s_{*R}(L) \]
hold obviously true.  Making use of the identities (\ref{lud}) and
(\ref{ftlint}) one  
checks that
\bea (\id_{\atld}\ot \gamma_{*R})\ci \gamma_{*L}(\phitl)&=&
\xi^{-2}(\ell\dii)\ru \phitl \ot \ell_{(1^{\prime})}\ru
{\ell_L}^{-1}(\ell\di)\ot 
{\ell_L}^{-1}(\ell_{(2^{\prime})})=\nn
&=&\xi^{-2}(\ell\dii)\ru \phitl \ot
{\ell_L}^{-1}\left(\ell\di\xi(\ell_{(1^{\prime})})\right)\ot 
{\ell_L}^{-1}(\ell_{(2^{\prime})})=\nn
&=&\xi^{-2}\left(\ell\dii\xi^2(\ell_{(1^{\prime})})\right)\ru \phitl
\ot {\ell_L}^{-1}(\ell\di)\ot 
{\ell_L}^{-1}(\ell_{(2^{\prime})})=\nn
&=&(\gamma_{*L}\ot \id_{\atld})\ci \gamma_{*R}(\phitl)\nn
(\gamma_{*R}\ot \id_{\atld})\ci\gamma_{*L}(\phitl)&=&
\ell_{(1^{\prime})}\xi^{-2}(\ell\dii)\ru \phitl \ot
{\ell_L}^{-1}(\ell_{(2^{\prime})})\ot 
{\ell_L}^{-1}(\ell\di)=\nn
&=&\ell_{(1^{\prime})}\ru \phitl \ot
{\ell_L}^{-1}\left(\ell_{(2^{\prime})}\xi^{-1}(\ell\dii)\right)\ot 
{\ell_L}^{-1}(\ell\di)=\nn
&=&\ell_{(1^{\prime})}\ru \phitl \ot \xi^{-2}(\ell\dii)\ru
{\ell_L}^{-1}(\ell_{(2^{\prime})})\ot 
{\ell_L}^{-1}(\ell\di)=\nn
&=& (\id_{\atld}\ot \gamma_{*L})\ci \gamma_{*R}(\phitl).
\nonumber\eea
The $S_*^{\ell}=({\ell_L}^{-1}\ci\xi \ci {_L\ell})\ci({_L\ell}^{-1}\ci
{\ell_L})$ is a composition of  
ring isomorphisms ${_L\ell}^{-1}\ci {\ell_L}:(\atld)^{op}\to \asld$
and ${\ell_L}^{-1}\ci \xi 
\ci {_L\ell}:\asld\to \atld$, hence it is an anti-automorphism of the
ring $\atld$. Also 
\bea S_*^{\ell}\ci t_{*R}(l)(a)&=&
\mu^{-1}\ci \ftr^{-1}\ci {\ell_L}\ci t_{*R}(l)\ci S^{-1}(a)=
\mu^{-1}\ci \ltr\left(S^{-1}[s_L(l
\pi_L(\ell\di))\ell\dii]S^{-1}(a)\right)\nn 
&=&\mu^{-1}\ci \pi_R\ci S^{-1}(a s_L(l))=\pi_L(as_L(l))=
s_{*R}(l)(a)\nn
S_*^{\ell}\ci t_{*L}(r)(a)&=&
\mu^{-1}\ci \ftr^{-1}\ci {\ell_L}\ci t_{*L}(r)\ci S^{-1}(a)=
\mu^{-1}\ci \lsr\left(a s_L\ci \pi_L(\ell\di t_R\ci
\kappa^{-1}(r))\ell\dii\right)\nn 
&=&\mu^{-1}\ci \lsr\left(t_L\ci \pi_L(a)\ell t_R\ci \kappa^{-1}(r)\right)=
\mu^{-1}(r)\pi_L(a)=s_{*L}(r)(a).
\nonumber\eea
Since 
\[ \gamma_L\ci \xi^{-1}(a)=
\xi^{-1}(a\uii)\ot S^{-1}(a\ui)\quad {\rm and}\quad 
\gamma_R\ci \xi^{-1}(a)= \xi^{-1}(a\dii)\ot S^{-1}(a\di)
\]
we have $\gamma_L\ci \xi^{-2}(a)=\xi^{-1}\ci S^{-1}(a\di)\ot S^{-1}\ci
\xi^{-1}(a\dii)$.
Then we can compute
\bea 
\left[S_*^{\ell-1}(\phitl\dii){\phitl}\di\right](a)=
\left[{\ell_L}^{-1}\ci
\xi^{-1}(\ell\di)(\xi^{-2}(\ell\dii)\ru\phitl)\right](a)= 
\!\!\!\!\!\!\!\!\!\!\!\!\!\!\!\!\!\!\!\!\!
\!\!\!\!\!\!\!\!\!\!\!\!\!\!\!\!\!\!\!\!\!
\!\!\!\!\!\!\!\!\!\!\!\!\!\!\!\!\!\!\!\!\!
\!\!\!\!\!\!\!\!\!\!\!\!\!\!\!\!\!\!\!\!\!
\!\!\!\!\!\!\!\!\!\!\!\!\!\!\!\!\!\!\!\!\!
\!\!\!\!\!\!\!\!\!\!\!\!\!\!\!\!\!\!\!\!\!
\!\!\!\!\!\!\!\!\!\!\!\!\!\!\!\!\!\!\!\!\!
&&\nn
&=&\phitl\left([a\lu {\ell_L}^{-1}\ci
\xi^{-1}(\ell\di)]\xi^{-2}(\ell\dii)\right)= \nn
&=&\phitl\left([\xi^{-1}(\ell\di)\ld
{_L\ell}^{-1}(a)]\xi^{-2}(\ell\dii)\right)= \nn
&=&\phitl\ci \xi^{-2}\left(\xi({\ell\uii}\di){\ell\uii}\dii\right)
{_L\ell}^{-1}(a) \ci S^{-1}(\ell\ui)=\nn
&=&\phitl(1_A) \pi_L(\ell\di){_L\ell}^{-1}(a)(\ell\dii)=
\phitl(1_A) \pi_L(a)
\nonumber\eea
hence $S_*^{\ell}({\phitl}\di){\phitl}\dii= s_{*R}\ci
\pi_{*R}(\phitl)$ for all $\phitl\in \atld$. Also
\bea [\phitl\ui S_*^{\ell}(\phitl\ui)](a)&=&
[(\ell\di \ru\phitl){\ell_L}^{-1}\ci \xi(\ell\dii)](a)=\nn
&=&\mu^{-1}\ci \ftr^{-1}(\ell\dii)\ci S^{-1}\left(s_L\ci \phitl(a\di
\ell\di) a\dii \right)=\nn
&=&\mu^{-1}\ci \lsr\left( s_L\ci \phitl(a\di\ell\di)a\dii\ell\dii\right)=
\mu^{-1}\ci \lsr\left(a\ell\lu\phitl\right)=\nn
&=&\mu^{-1}\ci \lsr\left(t_L\ci\pi_L(a)\ell\lu\phitl\right)=
\mu^{-1}\ci\lsr (\ell\lu\phitl)\pi_L(a)=
[s_{*L}\ci \pi_{*L}(\phitl)](a).
\nonumber\eea
This proves that ${\cal A}_*^{\ell}=({\cal
A}_{*L}^{\ell},\atld_R,S_*^{\ell})$ is a symmetrized Hopf algebroid. 
\hspace{1cm}\qed 

Obviously the strong isomorphism class of the Hopf algebroid $({\cal
A}_{*\ L}^{\ell},S_*^{\ell})$ depends only on the Hopf algebroid
$({\cal A}_L,S)$ and the non-degenerate left integral $\ell$ of it. It
is insensitive to the particular choice of the underlying right
bialgebroid ${\cal A}_R$.

The antipode $S_*^{\ell}$ has a form analogous to (\ref{xiarr}):
\be  S_*^{\ell}(\phitl)(a)=[(\ell\lu\phitl)\ru {\ell_L}^{-1}(1_A)](a).
\lb{S*arr}\ee
Using the left bialgebroid isomorphisms of Theorem \ref{bgdiso} also
the dual left bialgebroids $(\asld_R)^{op}_{cop}$, $\asrd_L$  and
$\atrd_L$ can be made 
Hopf algebroids all strictly isomorphic to the above Hopf algebroid
$({\cal A}_{* \ L}^{\ell},S_*^{\ell})$.  They have the antipodes
\bea {_*S^{\ell}}&=& {_L\ell}^{-1}\ci \xi\ci {_L\ell}\lb{Ssld}\\
     {S^*_{\ell}}&=&\fsr^{-1}\ci\xi\ci \fsr\lb{Ssrd}\\
     {^*S_{\ell}}&=&\ftr^{-1}\ci \xi\ci \ftr\lb{Strd}.
\eea
Let us turn to the interpretation of the role of the Hopf
algebroid $({\cal A}_{* L}^{\ell},S^{\ell}_*)$. As ${\cal A}_{*
L}^{\ell}$ is isomorphic to $(\asld_R)^{op}_{cop}$ and the right
bialgebroid underlying  $({\cal A}_{* L}^{\ell},S^{\ell}_*)$
is ${\cal A}_{* \ R}$,
on the first sight it seems to be natural to consider it
as some kind of a dual of $({\cal A}_L,S)$. There are however two
arguments 
against this interpretation: First, the Hopf algebroid  $({\cal A}_{*
L}^{\ell},S^{\ell}_*)$ 
depends on $\ell$, and it gives a generalization of the dual of a
finite dimensional Hopf algebra if and only if $S(\ell)=\ell$. Second,
as it is proven in the next Proposition \ref{d2snd},  $({\cal A}_{*
L}^{\ell},S^{\ell}_*)$  belongs to a
special kind of Hopf algebroids: it possesses a two sided
non-degenerate integral.
\bl \lb{finite}
Let ${\cal A}$ be a symmetrized Hopf algebroid such that the
$R$-module  $\asr$ in (\ref{rbim}) 
is finitely generated projective. Then a left integral $\ell\in \ila$
is non-degenerate if and only if the map $\fsr$ is bijective.
\el

\pr The {\em only if part} is trivial.

In order to prove the {\em if part} recall that by the proof of Lemma
\ref{lac} for $\lsr\colon =\fsr^{-1}(1_A)$ and all $a\in A$ the identity 
$ \lsr\ru a=s_R\ci \lsr (a)$ holds true, 
$\lsr$ is an $R$-$R$ bimodule map ${_R A^R}\to R$ and the inverse
of $\fsr$ reads as $\fsr^{-1}(a)=\lsr\lu S(a)$. Then
\[\phisr(a)=
\fsr^{-1}\ci \fsr(\phisr)(a)=
\lsr (s_R\ci \phisr(\ell\ui)S(\ell\uii)a)=
\phisr(a\ell\ui s_R\ci \lsr\ci S(\ell\uii))\]
for all $\phisr\in \asrd$ and $a\in A$. Using the finitely generated
projectivity of $\asr$, we have $\ell\ui s_R\ci \lsr\ci
S(\ell\uii)=1_A$. Since $\lsr\ci S\in \atrd$ the inverse $\ftr^{-1}$
can be defined as $\ftr^{-1}(a)=\lsr\ci S\ld S^{-1}(a)$.
\hspace{1cm}\qed

\bp \lb{d2snd}
Let ${\cal A}$ be a symmetrized Hopf algebroid possessing a non-degenerate left
integral $\ell$. Then the element ${\ell_L}^{-1}(1_A)$ in $\atld$ is
a two sided non-degenerate integral in the symmetrized Hopf algebroid ${\cal
A}_*^{\ell}$ (constructed in Theorem \ref{dualhgd}). 
\ep

\pr It follows from (\ref{ftl}) that ${\ell_L}^{-1}(1_A)=\mu^{-1}\ci
\lsr$ where $\lsr\colon =\fsr^{-1}(1_A)$ and $\mu$ is the ring
isomorphism introduced in Proposition \ref{sbgdi}.
For all $\phitl\in \atld$ and $a\in A$ we have
\bea [(\mu^{-1}\ci \lsr)\phitl](a)&=&
\phitl\ci S\left( \ltr\rd S^{-1}(a)\right)=
\phitl(1_A) \mu^{-1}\ci \lsr(a) \quad {\rm hence}\nn
\left[(\mu^{-1}\ci \lsr)t_{*R}\ci \pi_{*R}(\phitl)\right](a)&=&
t_{*R}\ci \pi_{*R}(\phitl)(1_A) \mu^{-1}\ci \lsr(a)=
\phitl(1_A)\mu^{-1}\ci \lsr(a)
\nonumber\eea
which proves that $\mu^{-1}\ci \lsr$ is a right integral. Using
(\ref{xiarr}) 
\[S_*^{\ell}(\mu^{-1}\ci\lsr)=
{\ell_L}^{-1}\ci \xi\ci {\ell_L}(\mu^{-1}\ci\lsr)=
{\ell_L}^{-1}\ci \xi^2(1_A)=
{\ell_L}^{-1}(1_A)=\mu^{-1}\ci\lsr\]
hence $\mu^{-1}\ci \lsr$ is also a left integral.

As it is proven in \cite{KSz}, since $\atl$ is finitely generated
projective so is the left $L_*\equiv L$-module ${^{L_*}(\atld)}$. 
The corresponding dual bialgebroid ${^*(\atld)}_L$ is isomorphic to
${\cal A}_L$ via the isomorphism $(\iota,\id_L)$ of left bialgebroids
where 
\be \iota:A\to {^*(\atld)}\qquad \iota(a)(\phitl)\colon = \phitl(a).
\lb{iota}\ee
Since
\bea (\iota(a)\rd\phitl)(b)&=&
\pi_L\left([b\lu(\ell\di\ru\phitl)]s_L\ci
{\ell_L}^{-1}(\ell\dii)(a)\right)=
\phitl\left(b s_R\ci \lsr(a \ell\ui)S(\ell\uii)\right)= \nn
&=&(a\ru\phitl)(b)
\nonumber\eea
the map  ${_{L_*}(\mu^{-1}\ci \lsr)}:{^*(\atld)}\to \atld$, $\iota(a)\
\mapsto\ \iota(a)\rd 
\mu^{-1}\ci \lsr \equiv a\ru \mu^{-1}\ci \lsr$ is bijective with
inverse
\[ {_{L_*}(\mu^{-1}\ci \lsr)}^{-1}:\phitl\ \mapsto\ \iota\ci {\ell_L}\ci
S_*^{\ell-1} \equiv \iota\ci \ftr(\mu\ci \phitl\ci S).\]
The application of Lemma \ref{finite} finishes the proof.
\hspace{1cm}\qed

In the view of  Proposition \ref{ndtw} the following definition
makes sense:

\bd \lb{hgddual}
The  {\em dual of the isomorphism class} of a  Hopf algebroid $({\cal
A}_L,S)$  possessing a non-degenerate left
integral $\ell$ is  the isomorphism class of the Hopf algebroid $({\cal
A}^{\ell}_{*L},S^{\ell}_*)$ (constructed in the Theorem \ref{dualhgd}). 
\ed
The next proposition shows that this notion of duality is involutive:

\bp Let  $\ell$  be  a non-degenerate left integral in  the Hopf
algebroid $({\cal A}_L,S)$. Then the Hopf algebroid 
$\left((\atld^{\ell})_{*\
L}^{{\ell_L}^{-1}(1_A)},(S_*^{\ell})_*^{{\ell_L}^{-1}(1_A)}\right)$  is
strictly isomorphic to $({\cal A}_L,S^{\prime}_{\ell})$ -- the Hopf
algebroid constructed in Proposition \ref{aprl}. In particular the
Hopf algebroid  $\left ((\atld^{\ell})_{*\
L}^{{\ell_L}^{-1}(1_A)},(S_*^{\ell})_*^{{\ell_L}^{-1}(1_A)}\right)$ is
isomorphic to $({\cal A}_L,S)$. 
\ep

\pr Since the Hopf algebroids $\left((\atld^{\ell})_{*\
L}^{{\ell_L}^{-1}(1_A)},(S_*^{\ell})_*^{{\ell_L}^{-1}(1_A)}\right)$
and  $({^*(\atld)}_L,{^*(S_*^{\ell})}_{{\ell_L}^{-1}(1_A)})$ are
strictly isomorphic it
suffices to show that the isomorphism $(\iota,\id_L)$ of left
bialgebroids ${\cal A}_L\to {^*(\atld)}_L$ in (\ref{iota}) extends to a
strict isomorphism of Hopf algebroids $({\cal A}_L,S^{\prime}_{\ell})
\to  ({^*(\atld)}_L,$ \break$ {^*(S_*^{\ell})}_{{\ell_L}^{-1}(1_A)})$.

By (\ref{Strd}) for a  non-degenerate left integral $\ell$ in the Hopf
algebroid $({\cal A}_L,S)$ the  antipode ${^* S_{\ell}}$ of the Hopf
algebroid $(\atrd_{L},{^* S_{\ell}})$ reads as 
\[{^* S_{\ell}}(\phitr)=\ltr\ld(\phitr\rd S^{-1}(\ell)).\]
Applying it to the non-degenerate integral ${{\ell_L}^{-1}(1_A)}$ in
$({\cal A}_{* L}^{\ell},S_*^{\ell})$  we obtain
\[ {^*(S_*^{\ell})}_{{\ell_L}^{-1}(1_A)}\left(\iota(a)\right)=
\iota(\ell)\ld \left[\iota(a)\rd S_*^{\ell-1}\ci {{\ell_L}^{-1}(1_A)}\right]=
\iota\left(\ell\lu (a\ru {{\ell_L}^{-1}(1_A)})\right)=
\iota\ci \xi(a). \qquad\qed \]

The duality of (weak) Hopf algebras is re-obtained from Definition
\ref{hgddual} as follows: Let $H$ be a
finite weak Hopf algebra  over a commutative ring $k$. Let $({\cal H}_L,S)$ 
be the corresponding Hopf algebroid --  introduced in the Example
\ref{whaex} -- , and let 
$\ell$ be a non-degenerate left integral in ${H}$.  Recall, that in order to
reconstruct the  weak Hopf algebra from the Hopf algebroid in Example
\ref{whaex} one needs a distinguished separability structure on the
base ring. The dual weak Hopf algebra is the unique weak Hopf algebra
in the isomorphism class of $({\cal H}_{* L}^{\ell},S_*^{\ell})$
corresponding to the same separability structure on $L$ as $H$
corresponds to.

If $H$ is a Hopf algebra over $k$ then -- since the separability
structure on $k$ is unique -- the dual Hopf algebra is the only Hopf
algebra in the isomorphism class of $({\cal H}_{*
L}^{\ell},S_*^{\ell})$. 


\begin{thebibliography}{NSzW}
\addcontentsline{toc}{section}{\protect\numberline{}{References}}

\bibitem{Barr} M. Barr: {\em $*$-autonomous categories}, Lecture Notes in
Mathematics {\bf 752}, Springer, Berlin, 1979

\bibitem{B} G. B\"ohm {\em `An alternative notion of Hopf algebroid'}
{\tt math.QA/0301169}

\bibitem{BNSz} G. B\"ohm, F. Nill, K. Szlach\'anyi {\em `Weak Hopf Algebras
I: Integral Theory and $C^*$-structure'} J. Algebra {\bf 221} (1999)
p. 385 

\bibitem{BSz} G. B\"ohm, K. Szlach\'anyi {\em `A Coassociative $C^*$-Quantum
Group with Non-Integral Dimensions'} Lett. Math. Phys. {\bf 35} (1996)
p. 137

\bibitem{BreMi} T. Brzezinski, G. Militaru {\em `Bialgebroids,
$\times_R$-bialgebras and Duality'} J. Algebra {\bf 247} No.2 (2002) p.467

\bibitem{ConnMo} A. Connes, H. Moscovici {\em `Cyclic cohomology and
Hopf algebra symmetry'} Conference Mosh Flato 1999, Dijon
Lett. Math. Phys. {\bf 52} No.1 (2000) p.1; {\em `Differential cyclic
cohomology and Hopf algebraic structures in transverse geometry'}
{\tt math.DG/0102167 }

\bibitem{DS} B. Day, R. Street {\em `Quantum categories, star
autonomy, and quantum groupoids'} \hfill\break{\tt math.CT/0301209}


\bibitem{KSz} L. Kadison, K. Szlach\'anyi {\em `Dual Bialgebroids for Depth
Two Ring Extensions'} \hfill\break{\tt math.RA/0108067}

\bibitem{K1} L. Kadison {\em `Hopf algebroids and H-separable
extensions'}  {\tt MPS 0201025}  to appear in Proc.\ Amer.\ Math.\ Soc.  

\bibitem{K2} L. Kadison {\em `A Hopf algebroid associated to a Galois
extension'} {\tt preprint}

\bibitem{KhalRang} M. Khalkhali, B. Rangipour {\em `On cohomology of
Hopf algebroids'} {\tt math.KT/0105105} to appear in Advances in
Mathematics 

\bibitem{LaSw} R.G. Larson, M.E. Sweedler {\em `An associative
orthogonal bilinear form for Hopf algebras'} Amer. J. of Math. {\bf
91} (1969) p. 75

\bibitem{Lu} J. H. Lu {\em `Hopf Algebroids and Quantum Groupoids'}
Int. J. Math. Vol. {\bf 7} No. {\bf 1} (1996) p. 47


\bibitem{N} F. Nill {\em `Weak Bialgebras'} {\tt math.QA/9805104}

\bibitem{Pa} F. Panaite {\em `Doubles of (quasi)Hopf algebras and some
examples of quantum groupoids and vertex groups related to them '}
{\tt math.QA/0101039 }


\bibitem{Schauenburg: Bial} P. Schauenburg {\em `Bialgebras over noncommutative
rings, and a structure theorem for Hopf bimodules'} Applied Categorical
Structures {\bf 6} (1998) p.193

\bibitem{Sch} P. Schauenburg {\em `Duals and Doubles of Quantum Groupoids'}
in: ``New trends in Hopf algebra theory'' Proceedings of the
Colloquium on quantum groups and Hopf algebras, La Falda, Sierra de
Cordoba, Argentina, 1999, AMS Contemporary Mathematics {\bf 267}
(2000) p. 273

\bibitem{Sch1} P. Schauenburg {\em `Weak Hopf Algebras and Quantum
Groupoids'} {\tt preprint}

\bibitem{Sw} M. E. Sweedler {\em `Groups of simple algebras'}
Publ. Math. I.H.E.S. {\bf 44} (1974) p.79

\bibitem{Sz} K. Szlach\'anyi: {\em `Finite Quantum Groupoids and Inclusions of
Finite Type'} Fields Institute Communications Vol {\bf 30} (2001) 393-407

\bibitem{Sz2} K. Szlach\'anyi {\em `Galois actions by finite quantum
groupoids'} in "Locally Compact Quantum Groups and Groupoids", proceedings of
the meeting of theoretical physicists and mathematicians, Strasbourg, February
2002,  ed.: L. Vainerman  (IRMA Lectures in Mathematics and Theoretical
Physics 2, series editor: V. Turaev) de Gruyter 2003

\bibitem{Sz: EM} K. Szlach\'anyi {\em The monoidal Eilenberg-Moore
construction and bialgebroids}, to appear in J. Pure Appl. Algebra

\bibitem{T} M. Takeuchi {\em `Groups of Algebras over $A\otimes {\bar A}$'}
J. Math. Soc Japan {\bf 29} (1997) p. 459

\bibitem{Xu} P. Xu {\em `Quantum Groupoids and Deformation Quantization'}
C.R. Acad. Sci. Paris, I. {\bf 326} (1998) p. 289 {\em `Quantum
Groupoids'} {\tt math.QA/9905192}

\bibitem{V} P. Vecserny\'es {\em `Larson-Sweedler theorem, group like
elements, invertible modules and the order of the antipode in weak Hopf algebras'}
{\tt math.QA/0111045}

\bibitem{W} Y. Watatani {\em `Index for $C^*$-subalgebras'} Memoirs of
the AMS Vol. 23 No. 424 (1990)

\end{thebibliography}
\end{document}